\documentclass[preprint,12pt,authoryear]{elsarticle}

\usepackage{epstopdf}
\usepackage{graphicx}%
\usepackage{multirow}%
\usepackage{amsmath,amssymb,amsfonts}%
\usepackage{amsthm}%
\usepackage{mathrsfs}%
\usepackage[title]{appendix}%
\usepackage{xcolor}%
\usepackage{textcomp}%
\usepackage{manyfoot}%
\usepackage{booktabs}%
\usepackage{algorithm}%
\usepackage{algorithmic}%
\usepackage{listings}%
\usepackage{adjustbox}
\usepackage{subfig}
\usepackage{pifont}
\usepackage{colortbl}
\usepackage[a4paper,margin=25mm]{geometry}
\usepackage[hidelinks,hypertexnames=false]{hyperref}
\biboptions{round,semicolon}
\hypersetup{pdfauthor={Zhuoxin Fan, Liping Tang, Xinmin Yang},pdftitle={An Inexact Proximal Gradient Method for Composite Multiobjective Optimization}}

\theoremstyle{plain}
\newtheorem{theorem}{Theorem}
\newtheorem{lemma}{Lemma}
\newtheorem{assumption}{Assumption}
\newtheorem{corollary}{Corollary}
\newtheorem{proposition}[theorem]{Proposition}

\theoremstyle{definition}

\newtheorem{remark}{Remark}
\newtheorem{definition}{Definition}

\newcommand{\R}{\mathbb{R}}

\newcommand{\norm}[1]{\left\lVert #1 \right\rVert}

\newcommand{\inner}[2]{\left\langle #1,\,#2 \right\rangle}

\newcommand{\prox}{\operatorname{prox}}
\newcommand{\dom}{\text{dom}}

\journal{European Journal of Operational Research}

\begin{document}
\begin{frontmatter}

\title{An Inexact Proximal Gradient Method for Composite Multiobjective Optimization with Absolute and Relative Stopping Criteria}

\author[aff1]{Zhuoxin Fan}
\ead{fanzx@163.com}

\author[aff1]{Liping Tang}
\ead{tanglipings@163.com}

\author[aff1]{Xinmin Yang\corref{cor1}}
\ead{xmyang@cqnu.edu.cn}
\cortext[cor1]{Corresponding author.}

\affiliation[aff1]{organization={National Center for Applied Mathematics in Chongqing, Chongqing Normal University},
                  city={Chongqing},
                  postcode={401331},
                  country={China}}

\begin{abstract}
Proximal gradient methods are widely used for composite optimization. In multiobjective optimization, however, the resulting minimax subproblem generally requires an additional inner method to solve, even when the proximal mapping of each nonsmooth component is available in closed form. In addition, most existing convergence results are obtained only under the assumption that the involved subproblems are solved exactly, which is unrealistic and limits the applicability of the proximal gradient methods. To make the proximal gradient methods practical, we propose an implementable inexact proximal gradient method based on proximal computation. Considering an approximate optimality condition and a hybrid stopping criterion including absolute and relative rules for the subproblem, we address the computational challenge by developing an inner inexact projected gradient method that can exploit the proximal mappings of the individual components. The $\mathcal{O}(1/k)$ convergence rate is also established for our method under certain conditions. Numerical experiments illustrate the computational benefits of solving the subproblems inexactly.
\end{abstract}


\end{frontmatter}

\section{Introduction}
In this paper, we consider the following composite multiobjective optimization problem:
\begin{equation}\label{def:MOP}
	\min_{x\in\R^n} F(x)=f(x)+g(x)+h(x)\boldsymbol 1,
	\tag{MOP}
\end{equation}
where $f=(f_1,\ldots,f_m):\R^n\to\R^m$ and $g=(g_1,\ldots,g_m):\R^n\to\R^m$. Each $f_i$ is convex, continuously differentiable and $L$-smooth, each $g_i$ is convex and $L_i$-Lipschitz continuous on $\R^n$, and $h:\R^n\to\R\cup\{+\infty\}$ is proper, lower semicontinuous and convex. Here $\boldsymbol 1=(1,\ldots,1)\in\R^m$, so $F_i=f_i+g_i+h$ and $\dom F=\dom h\neq\emptyset$. The common term $h$ may include the indicator function of a nonempty closed convex feasible set.

Composite multiobjective optimization problems arise naturally in applications where the objective functions contain nonsmooth regularizers. Each component of the vector-valued objective may involve a differentiable convex term together with a nonsmooth convex term, such as an indicator function of a feasible set or a regularizer promoting certain structures. Such problems appear in robust multiobjective optimization \citep{TanabeFukudaYamashita2019}.

Two classical approach are scalarization and metaheuristics method. The scalarization includes the weighted-sum method \citep{Miettinen1999}, the $\epsilon$-constraint method \citep{Miettinen1999}, and the Pascoletti--Serafini scalarization \citep{PascolettiSerafini1984}. Scalarization is theoretically important and allows the use of scalar optimization techniques. However, the quality of the obtained solutions usually depends on the choice of scalarization parameters, which may be difficult to determine in advance. Another class of approaches is based on metaheuristics, such as NSGA-II \citep{Deb2002} and MOEA/D \citep{ZhangLi2007}. Although these methods are widely used in practice, they generally do not provide the same type of deterministic convergence guarantees as descent-type methods.

Recently, descent methods form an important class of direct methods for multiobjective optimization. In smooth multiobjective optimization, Fliege and Svaiter \citeyearpar{FliegeSvaiter2000} proposed a steepest descent method that computes a common descent direction for all objective functions at each iteration. This idea has been extended to conjugate gradient methods \citep{LucambioPerezPrudente2018}, Newton methods \citep{FliegeDrummondSvaiter2009}, and projected gradient methods \citep{DrummondIusem2004}. Further developments include nonmonotone gradient methods for vector optimization \citep{QuJiJiangZhang2017} and a Barzilai--Borwein descent method for multiobjective optimization \citep{MorovatiPourkarimiBasirzadeh2016,ChenTangYang2023}. For composite multiobjective optimization, proximal gradient methods have also been developed to handle the coexistence of smooth and nonsmooth terms \citep{TanabeFukudaYamashita2019}.

Although proximal gradient methods are effective for composite multiobjective optimization, their practical implementation requires solving a subproblem at each iteration. Define $\phi_i(u;x):=\inner{\nabla f_i(x)}{u-x}+g_i(u)-g_i(x)$ and
$\phi(u;x):=(\phi_1(u;x),...,\phi_m(u;x))$. The common term $h$ is kept outside the maximum. For $x\in\dom F$ and $\beta>0$, the subproblem \citep{TanabeFukudaYamashita2019} is
\begin{equation}\label{eqn:primal}
	\min_{u\in\R^n}\Bigl\{
	\Phi_{p}(u;x):= \max_{i\in [m]}\phi_i(u;x)+h(u)-h(x)+\frac{1}{2\beta}\norm{u-x}^2
	\Bigr\}.
	\tag{P}
\end{equation}
In general, an exact solution of \eqref{eqn:primal} is not available in closed form and must be computed by an iterative method. In particular, even when the proximal mappings of $g_i$ and $h$ admit closed-form expressions, this does not generally yield a closed form solution of the subproblem, since the proximal mapping associated with $h+\sum_i\lambda_i g_i$ cannot, in general, be recovered directly from their individual proximal mappings. Therefore, this leads us to develop a implementable and practical proximal gradient method for multiobjective optimization, and establish the convergence analysis under inexactness \citep{Schmidt2011,BelloCruz2025Relative,FukudaDrummond2013}.

The main contributions of this paper can be summarized as follows.

(1) We introduce an approximate optimality condition (AOC) and a stopping criterion for the proximal gradient subproblem in composite multiobjective optimization. The proposed AOC provides a unified characterization of inexact subproblem solutions by simultaneously accounting for errors in the proximal computations and in the descent direction. Based on the AOC, we further propose an stopping criterion for controlling the inexactness and an Armijo type line search that accommodates such errors.

(2) We develop an inexact projected gradient method for the dual form of \eqref{eqn:primal} that returns a subproblem solution satisfying both the AOC and the inexact stopping criterion. We can obtain the inexact gradient of the dual subproblem by using ADMM to approximate the proximal point for $h+\sum_i\lambda_i g_i$ using the individual proximal operators of $h$ and the $g_i$. Hence, our method is based on proximal computation.

(3) We establish convergence results for the proposed methods. In particular, we analyze the convergence of the inner inexact projected gradient method and show that it can generate approximate subproblem solutions satisfying the required conditions. We then establish the convergence of the outer inexact proximal gradient method under the proposed stopping criterion. Numerical experiments are also presented to illustrate the practical performance of the proposed methods. 

The paper is organized as follows. In Sect. \ref{sec:pre} we introduce the notation and some lemmas that will be used in the remainder of the paper. In Sect. \ref{sec:inexact} we present an inexactness analysis of the subproblem and propose AOC and a stopping criterion. In Sect. \ref{sec:pg} we introduce an explicit linesearch and the inexact proximal gradient method for multiobjective optimization. The inexact projected gradient method also be presented in this section. The convergence analyses of these two algorithms are presented in Sect. \ref{sec:conv}. Section 6 presents numerical examples and comparisons. In Sect. \ref{sec:conclusion} we present some conclusions.

\section{Preliminaries}
\label{sec:pre}
In this section, we present some preliminary material and notations that will be used throughout this paper. Denote $[m]=\{1,2,\ldots,m\}$, $\R_+^n:=\{x\in\R^n:x\geq 0\}$, $\R_{++}^n:=\{x\in\R^n:x > 0\}$ and $\Delta_m:=\{\lambda\in \R^m_+:\sum_{i=1}^{m}\lambda_i=1\}$. The definition of the $\epsilon$-normal cone of $x\in C\subseteq \R^n$ is
\[
N_C^\epsilon(x)
=
\left\{
v\in\mathbb{R}^n :
\langle v, y-x\rangle \le \epsilon,
\quad \forall y\in C
\right\},
\quad \epsilon \ge 0.
\]
Throughout this paper, unless otherwise specified, the norm $\norm{\cdot}$ is the Euclidean norm. Let $f:\R^n\to\R\cup\{+\infty\}$ be a proper, lower semicontinuous (l.s.c.), and convex function. We give some classical results of $f(x)$. At first, we give the definition of $\epsilon$-subdifferential.
\label{sec:preliminaries}
\begin{definition}[$\epsilon$-subdifferential]
	The $\epsilon$-subdifferential of the proper l.s.c function $f:\R^n\to \R\cup\{+\infty\}$ at the point $x\in\dom f$ is the set
	\begin{equation}
		\partial_\epsilon f(x)=\{\xi\in\R^n:f(y)\geq f(x)+\inner{y-x}{\xi}-\epsilon,\forall y\in\dom f\}.
	\end{equation}
\end{definition}
It is worth to note that
\begin{align*}
	w\in\partial_\epsilon f(x)\iff f(x)&\leq \inf_{y\in\mathcal{C}}\left\{f(y)+\inner{x-y}{w}\right\}+\epsilon\\
	&\leq\inf_{y\in\mathcal{C}}\left\{f(y)+\frac{1}{2}\norm{x-y}^2\right\}+\frac{1}{2}\norm{w}^2+\epsilon.
\end{align*}
When $w=0$, we get the equivalent relation between $\epsilon$ approximate solution and $\epsilon$ approximate optimality condition,
\[0\in\partial_\epsilon f(x)\iff f(x)\leq \inf_{y\in\dom f}f(y)+\epsilon.\]
For a proper, upper semicontinuous, concave function $\psi$, we can define the $\epsilon$-superdifferential.
\begin{definition}[$\epsilon$-superdifferential]
	Let $\psi:\R^n\to\R\cup\{-\infty\}$ be a proper,
	upper semicontinuous, concave function. For $x\in\dom \psi$
	and $\epsilon\geq0$, define
	\[
	\partial_\epsilon^+\psi(x)
	=
	\left\{
	\xi\in\R^n:
	\psi(y)\leq \psi(x)+\langle\xi,y-x\rangle+\epsilon,
	\quad \forall y\in\dom \psi
	\right\}.
	\]
\end{definition}
The following transportation formula will be used to express
the componentwise proximal optimality conditions at a common
point in the ADMM procedure; see
\citep[Chapter~XI]{HiriartUrrutyLemarechal1993II}.
\begin{proposition}
	\label{prop:subdiff}
	Let $v\in\partial f(y)$ where $x,y\in \dom f$. Then $v\in\partial_{\epsilon}f(x)$ where
	$$\epsilon=f(x)-f(y)-\inner{v}{x-y}\geq 0.$$
\end{proposition}
The following proximal three-point inequality will be used
to estimate the error in approximate proximal computations;
see also \citep{ChenTeboulle1993}.
\begin{lemma}
	Suppose that $v\in\prox_f(z)$, then we have $z-v\in\partial f(v)$ and for all $y\in\dom f$, we have
	\begin{equation}
		f(y)+\frac{1}{2}\norm{y-z}^2-\left(f(v)+\frac{1}{2}\norm{v-z}^2\right)\geq\frac{1}{2}\norm{y-v}^2.
	\end{equation}
	If $z-v\in\partial_{\epsilon} f(v)$, then for all $y\in\dom f$, we have
	\begin{equation}
		f(y)+\frac{1}{2}\norm{y-z}^2-\left(f(v)+\frac{1}{2}\norm{v-z}^2\right)\geq\frac{1}{2}\norm{y-v}^2-\epsilon.
	\end{equation}
\end{lemma}
\begin{proof}
	Suppose that $z-v\in\partial_\epsilon f(v)$. Then, for every
	$y\in\dom f$,
	\[
	f(y)-f(v)\geq \langle z-v,y-v\rangle-\epsilon.
	\]
	Adding the identity
	$\frac12\|y-z\|^2-\frac12\|v-z\|^2=\langle v-z,y-v\rangle+\frac12\|y-v\|^2$
	gives the second inequality. The first follows by taking
	$\epsilon=0$ and using
	$v=\prox_f(z)\iff z-v\in\partial f(v)$.
\end{proof}
\section{Inexactness analysis}
\label{sec:inexact}
In this section, we analyze the inexactness of the subproblem for the proximal gradient method. We first study its first order optimality conditions and then introduce an approximate optimality condition. We briefly show that some existing stopping criteria are covered by our approximate optimality condition. Finally, based on this condition, we propose a stopping criterion.

The associated saddle point problem of \eqref{eqn:primal} is
\begin{equation}\label{eqn:saddle}
	\max_{\lambda\in\Delta_m}\min_{u\in\R^n}\Bigl\{
	\Phi(u,\lambda;x):= \sum_{i=1}^{m}\lambda_i\phi_i(u;x)+h(u)-h(x)+\frac{1}{2\beta}\norm{u-x}^2
	\Bigr\}.
	\tag{SP}
\end{equation}
The epigraph formulation of~\eqref{eqn:primal} has a strictly feasible point with $u\in\operatorname{ri}(\dom h)$, since each $\phi_i$ is finite and continuous. Hence strong duality holds. We get the associated dual problem
\begin{equation}
	\label{eq:dual-1}
	\max_{\lambda\in\Delta_m}\Bigl\{
	\Phi_d(\lambda;x):=\min_{u\in\R^n}\Phi(u,\lambda;x)
	\Bigr\}.
	\tag{D}
\end{equation}
Denote that $\hat{x}=x-{\beta}\sum_{i=1}^{m}\lambda_i\nabla f_i(x)$ and $\mathcal{M}_\psi(x)=\min_{y\in\R^n}\left\{\psi(y)+\frac{1}{2}\norm{x-y}^2\right\}$ for a proper l.s.c. convex function $\psi$. We have the equivalent form
\begin{equation}\label{eq:dual-2}
	\Phi_{d}(\lambda;x)=\frac{1}{\beta}\mathcal{M}_{\beta(h+\sum_{i=1}^{m}\lambda_i g_i)}(\hat{x})-\frac{\beta}{2}\norm{\sum_{i=1}^{m}\lambda_i\nabla f_i(x)}^2-\sum_{i=1}^{m}\lambda_ig_i(x)-h(x).
\end{equation}

Fix $x\in\dom F$ and $\beta>0$. Whenever a saddle point $(u^\star,\lambda^\star)$ exists, its optimality conditions are as follows. Let $(u^\star,\lambda^\star)$ be the optimal solution of \eqref{eqn:saddle}. The first order optimality condition is:
\begin{equation*}
	\left\{
	\begin{aligned}
		0 &\in u^\star - x + \beta\sum_{i=1}^m \lambda_i^\star\nabla f_i(x) +\beta \partial \Bigl(h+\sum_{i=1}^m \lambda_i^\star g_i\Bigr)(u^\star),\\
		0 &\in -\phi(u^\star;x) + N_{\Delta_m}\bigl(\lambda^\star\bigr).
	\end{aligned}
	\right.
\end{equation*}
The first inclusion characterizes the proximal solution
for fixed weights, whereas the second expresses optimality
with respect to the weights on the simplex.
We discuss the inexactness of the proximal operator. For a proper l.s.c. convex function $\psi$ and $u\in\prox_\psi(z)$, we have $z-u\in\partial\psi(u)$.
If the proximal operator cannot be solved exactly, we may instead obtain
\begin{equation}
	\begin{cases}
		w\in\partial_\epsilon\psi(u)-z + u, \\
		\norm{w} \text{ and } \epsilon \text{ are small}.
	\end{cases}
\end{equation}
We relax the condition on the weights through an
$\epsilon_2$-normal cone. This motivates the following approximate optimality condition. 
\begin{definition}
	We say that $(\tilde{u},\tilde{w},\tilde{\lambda},\epsilon)\in\R^n\times\R^n\times\Delta_m\times\R^2_+$ satisfies the approximate optimality condition \ref{def:AOC} at $x\in\dom F$ if $\tilde u,\tilde u-\beta\tilde w\in\dom h$ and the following relation holds,
	\begin{equation}\label{def:AOC}
		\left\{
		\begin{aligned}
			\tilde{w} &\in \frac{1}{\beta}(\tilde{u} - x) + \sum_{i=1}^m \tilde{\lambda}_i\nabla f_i(x) + \partial_{\epsilon_1} \Bigl(h+\sum_{i=1}^m \tilde\lambda_i g_i\Bigr)(\tilde u),\\
			0 &\in -\phi(\tilde{u}-\beta \tilde{w};x) + N_{\Delta_m}^{\epsilon_2}\bigl(\tilde{\lambda}\bigr).
		\end{aligned}
		\right.
		\tag{AOC}
	\end{equation}
\end{definition}

\begin{remark}
	The AOC covers existing inexact conditions in both
	multiobjective and single-objective optimization.
	When $m>1$, $g_i=0$ for all $i\in[m]$, and $h=\delta_C$,
	the sufficient condition of
	\citet{FukudaDrummond2013} fits the AOC by taking
	$w=0$, $\epsilon_1=0$, and
	$\epsilon_2=-\sigma\Phi_d(\lambda;x)$,
	where $\sigma\in[0,1)$ is their relative error parameter.
	When $m=1$, the second inclusion is automatic with
	$\epsilon_2=0$, and the first reduces to
	\[
	w\in\nabla f(x)+\partial_{\epsilon_1}(g+h)(u)
	+\frac{1}{\beta}(u-x),
	\]
	which, for $\beta=1$ and nonsmooth term $g+h$, coincides with the approximate
	optimality inclusion used by
	\citet{BelloCruz2025Relative}.
\end{remark}
The following lemma states that if the first inclusion relation of \ref{def:AOC} is satisfied, then an important inequality holds.
\begin{lemma}\label{lem:pgap}
	Suppose that
	\begin{equation*}
		w\in\sum_{i=1}^{m}\lambda_i\left(\nabla f_i(x)+\frac{1}{\beta}(u-x)\right)+\partial_{\epsilon_1}\left(h+\sum_{i=1}^{m}\lambda_i g_i\right)(u),
	\end{equation*}
	Then for all $z\in\dom h$, we have
	\begin{align}\label{eqn:firstbd}
			&h(u)-h(z)+\sum_{i=1}^{m}\lambda_i\left(g_i(u)-g_i(z)+\inner{\nabla f_i(x)}{u-z}\right)\nonumber\\
			&\leq\frac{1}{2\beta}\left(\norm{z-x}^2-\norm{u-x}^2-\norm{u-z}^2\right)
			+\inner{w}{u-z}+\epsilon_1.
	\end{align}
	Equivalently, for all $z\in\dom h$, we have
	\begin{equation}\label{eqn:pdgap}
		\Phi(u,\lambda;x)-\Phi(z,\lambda;x)\leq-\frac{1}{2\beta}\norm{u-z}^2
		+\inner{w}{u-z}+\epsilon_1.
	\end{equation}
\end{lemma}
\begin{proof}
    From the definition of $\epsilon$-subdifferential, for all $z\in\dom h$,
    \[
    h(u)-h(z)+\sum_{i=1}^{m}\lambda_i\left(g_i(u)-g_i(z)+\inner{\nabla f_i(x)}{u-z}\right)
    \leq\inner{w}{u-z}-\frac1\beta\inner{u-x}{u-z}+\epsilon_1.
    \]
    Using $2\inner{u-x}{u-z}=\norm{u-x}^2+\norm{u-z}^2-\norm{z-x}^2$
    gives~\eqref{eqn:firstbd}. Adding the difference of the quadratic terms in $\Phi$
    proves~\eqref{eqn:pdgap}.
\end{proof}
From Lemma \ref{lem:pgap}, we can get the relation between the duality gap and the errors.
Denote that $\hat{x}=x - \beta \sum_{i=1}^{m} \lambda_i \nabla f_i(x)$, and $v=\prox_{\beta(h+\sum_{i=1}^{m}\lambda_i g_i)}\left(\hat{x}\right)$.
The duality gap of \eqref{eqn:saddle} is
\begin{align*}
		&\Phi_{p}(u-\beta w;x)-\Phi_{d}(\lambda;x)\\
		=&\max_{i\in [m]}\phi_i(u-\beta w;x)-\sum_{i=1}^{m}\lambda_i\phi_i(v;x)+h(u-\beta w)-h(v)+\frac{1}{2\beta}\norm{u-\beta w-x}^2-\frac{1}{2\beta}\norm{v-x}^2\\
		=&\max_{i\in [m]}\phi_i(u-\beta w;x)-\sum_{i=1}^{m}\lambda_i\phi_i(u-\beta w;x)+\Phi(u,\lambda;x)-\Phi(v,\lambda;x)\\
		&+\Phi(u-\beta w,\lambda;x)-\Phi(u,\lambda;x).
\end{align*}
By \eqref{eqn:pdgap} and Young's inequality, we have
\begin{equation*}
		\Phi(u,\lambda;x)-\Phi(v,\lambda;x)\leq-\frac{1}{2\beta}\norm{u-v}^2
		+\inner{ w}{u-v}+\epsilon_1
		\leq\frac{\beta}{2}\norm{w}^2+\epsilon_1.
\end{equation*}
Therefore, by \ref{def:AOC}, the following corollary holds.

\begin{corollary}
	Suppose that $(u,w,\lambda,\epsilon)$ satisfies \ref{def:AOC} and the assumptions of Lemma~\ref{lem:pgap}. Then the duality gap of \eqref{eqn:saddle} satisfies
	\begin{equation}\label{eqn:pdgap-ubd}
		\Phi_{p}(u-\beta w;x)-\Phi_{d}(\lambda;x)\leq\frac{\beta}{2}\norm{w}^2+\Phi(u-\beta w,\lambda;x)-\Phi(u,\lambda;x)+\epsilon_1+\epsilon_2.
	\end{equation}
\end{corollary}
The preceding corollary bounds the duality gap in terms of
the approximation errors.
However, \ref{def:AOC} only describes approximate optimality
for given residuals and tolerances; it does not prescribe how
small they must be.
Without further control, these errors may offset the decrease
obtained at each iteration and prevent convergence.
This motivates a stopping criterion to use in inexact proximal methods
\citep{Schmidt2011,BelloCruz2025Relative}.
We therefore bound the error by a relative term and a summable tolerance.
This ties the required accuracy to the progress of the outer
iteration and limits the effect of the remaining errors.
\begin{definition}[Inexact solution]
	Suppose that $(u,w,\lambda,\epsilon)$ satisfies \ref{def:AOC}. Define the inexact residual  $\xi(x)$ at $x$ as follows:
	\begin{equation}
		\begin{aligned}
			\xi(x)=&\frac{\tau}{2\beta}\norm{u-\beta w-x}^2-\frac{\tau}{2\beta}\norm{x-u}^2+\frac{\beta}{2}\norm{w}^2\\
			&+\sum_{i=1}^{m}\lambda_i\left(g_i(u-\beta w)-g_i(u)-\beta \inner{\nabla f_i(x)}{w}\right)\\
            &+h(u-\beta w)-h(u)+2\epsilon_1+\epsilon_2.
		\end{aligned}
	\end{equation}
	For $k\geq 0$, we say $u^k$ is an inexact solution at $x^k\in\dom F$ if there exist $(u,w,\lambda,\epsilon)\in\R^n\times\R^n\times\R^m\times\R_+^2$, $0<\tau<1$ and $0\leq\tau_0<1$ such that $(u^k,w^k,\lambda^k,\epsilon^k)$ satisfies \ref{def:AOC} and
	\begin{equation}\label{eq:xi-ubd}
		\xi^k\leq \frac{(1-\tau)\tau_0}{2\beta_k}\norm{u^k-x^k}^2+\eta_k,
	\end{equation}
	where $\{\eta_k\}$ is a nonnegative, monotone nonincreasing sequence and $\sum_{k=0}^{+\infty}\eta_k<+\infty$.
\end{definition}
\begin{remark}
	The bound on $\xi^k$ combines relative and absolute error control. When $m=1$ and $\eta_k\equiv0$, condition~\eqref{eq:xi-ubd} yields a relative stopping criterion similar to that used in the definition of inexact proximal solutions by \citet{BelloCruz2025Relative}. When $m=1$ and $\tau_0=0$, the relative term vanishes and the condition becomes $\xi^k\leq\eta_k$, an absolute stopping criterion with a prescribed summable tolerance. Related absolute stopping criteria for inexact proximal gradient methods have been studied by \citet{Schmidt2011}.
\end{remark}
We can also derive a lower bound for $\xi(x)$ under Lemma \ref{lem:pgap}.
\begin{corollary}
	Suppose that the assumptions of Lemma~\ref{lem:pgap} hold. $\xi(x)$ is the inexact residual at $x\in\dom F$. Then it satisfies the following lower bound,
	\begin{equation}\label{eq:reslbd}
		\xi(x)\geq \epsilon_1+\epsilon_2+\frac{\tau-1}{2\beta}\norm{u-\beta w-x}^2+\frac{1-\tau}{2\beta}\norm{u-x}^2.
	\end{equation}
\end{corollary}
\begin{proof}
	Setting $z=u-\beta w$ in the inequality \eqref{eqn:firstbd}, we get
	\begin{align*}
		&h(u)-h(u-\beta w)+\sum_{i=1}^{m}\lambda_i\left(g_i(u)-g_i(u-\beta w)+\inner{\nabla f_i(x)}{\beta w}\right)\\
		\leq&\frac{1}{2\beta}\left(\norm{u-\beta w-x}^2-\norm{u-x}^2-\norm{\beta w}^2\right)
		+\inner{w}{\beta w}+\epsilon_1\\
		=&\frac{1}{2\beta}\norm{u-\beta w-x}^2-\frac{1}{2\beta}\norm{u-x}^2+\frac{\beta}{2}\norm{w}^2+\epsilon_1.
	\end{align*}
	Then from the definition of $\xi(x)$, we have
	\begin{equation*}
		\xi(x)\geq \epsilon_1+\epsilon_2+\frac{\tau-1}{2\beta}\norm{u-\beta w-x}^2+\frac{1-\tau}{2\beta}\norm{u-x}^2.
	\end{equation*}
\end{proof}
\section{The inexact proximal gradient method}
\label{sec:pg}
This section presents the proposed method and its inner solvers. The outer proximal gradient method requires an approximate solution satisfying the \ref{def:AOC} and the stopping criterion. An inexact projected gradient method
generates this solution by solving the dual subproblem. The required weighted proximal subproblem of the inner solver are then solved by ADMM using the individual proximal operators.
\subsection{Inexact proximal gradient method for MOP}
The linesearch and proposed method are formally described below.
\begin{algorithm}[htbp]
	\begin{normalsize}
		\caption{Inexact proximal gradient method for
			multiobjective optimization}
		\label{algo:proximal}
		\begin{algorithmic}
			\REQUIRE Initial estimate $x^0\in\dom F$,
			parameter $\beta_0>0$, the optimization method
			$\mathcal P$, and tolerance $\bar\epsilon>0$.
			
			\STATE Set $k=0$.
			
			\WHILE{true}
			\STATE Find
			$(u^k,w^k,\lambda^k,\epsilon^k)
			\in\mathcal P(x^k,\beta_k)$
			such that \ref{def:AOC} and \eqref{eq:xi-ubd} hold;
			
			\IF{$\norm{u^k-\beta_kw^k-x^k}<\bar\epsilon$}
			\RETURN $x^k$;
			\ENDIF
			
			\STATE Find $\gamma_k\in(0,1]$ with line search;
			
			\STATE Update
			$x^{k+1}=x^k+\gamma_k(u^k-\beta_kw^k-x^k)$;
			
			\STATE Choose $\beta_{k+1}\in\R_{++}$;
			
			\STATE $k\leftarrow k+1$;
			\ENDWHILE
		\end{algorithmic}
	\end{normalsize}
\end{algorithm}

\paragraph{\textbf{Line search}}
Define
\begin{align}
	\label{eq:ls-omega}	\Omega^k=&\sum_{i=1}^{m}\lambda^k_i\left(g_i(u^k-\beta_kw^k)-g_i(x^k)\right)+\inner{\sum_{i=1}^{m}\lambda^k_i\nabla f_i(x^k)}{u^k-\beta_kw^k-x^k}\nonumber\\
	&+h(u^k-\beta_kw^k)-h(x^k)+\frac{\tau}{2\beta_k}\norm{u^k-\beta_kw^k-x^k}^2+\epsilon^k_1+\epsilon_2^k,
\end{align}
where $0<\tau<1$. Find the largest $\gamma=\frac{1}{2^j},j=0,1,...$ such that the following inequality holds:
\begin{equation*}
	\max_{i\in [m]}\left(F_i(x^k + \gamma(u^k -\beta_kw^k - x^k))-F_i(x^k)\right)\leq\gamma \Omega^k.
\end{equation*}
\begin{remark}
	This line search follows the fixed direction backtracking strategy of the Type~2 line search (LS2) in \citet{Salzo2017} when $w=0$ and $\epsilon_1=\epsilon_2=0$. Related inexact explicit line searches were studied in \citet{BelloCruz2025Relative,BelloCruz2025MOP}.
\end{remark}
Now, we state that the line search can stop after finite iterations, and we can get the lower bound of $\gamma_k$.
\begin{lemma}
	Let $0<\beta_{\min}\leq\beta_k\leq\beta_{\max}<+\infty$ for all $k\geq 0$. Define $d^k=u^k-\beta_k w^k-x^k$ and $\Omega^k$ as \eqref{eq:ls-omega}. The line search is well defined. Moreover, if $\tau>0$, then the step sizes admit the uniform lower bound
	$\gamma_k\ge\gamma_{\min}:=\min\left\{1,\frac{\tau}{2\beta_{\max}L}\right\}$, for all $k\ge 0$.
	
\end{lemma}

\begin{proof}
	Fix any $k\ge 0$. Since each $g_i$ is convex, for any $\gamma\in(0,1]$, we have $g_i(x^k+\gamma d^k)\le (1-\gamma)g_i(x^k)+\gamma g_i(x^k+d^k)$,
	which implies $g_i(x^k+\gamma d^k)-g_i(x^k)	\le\gamma\bigl(g_i(x^k+d^k)-g_i(x^k)\bigr)$. Since $h$ is convex and $x^k,x^k+d^k\in\dom h$, the same inequality holds with $g_i$ replaced by $h$. All trial points therefore belong to $\dom h$.
	On the other hand, since $\nabla f_i$ is $L$-Lipschitz continuous,
	\[
	f_i(x^k+\gamma d^k)-f_i(x^k)
	\le \gamma\inner{\nabla f_i(x^k)}{d^k}
	+\frac{L}{2}\gamma^2\|d^k\|^2.
	\]
	Therefore,
	\begin{align*}
	F_i(x^k+\gamma d^k)-F_i(x^k)
	&\le \gamma\Bigl(g_i(x^k+d^k)-g_i(x^k)+\inner{\nabla f_i(x^k)}{d^k}\Bigr)\\
	&\quad+\gamma\bigl(h(x^k+d^k)-h(x^k)\bigr)+\frac{L}{2}\gamma^2\|d^k\|^2.
	\end{align*}
	Taking the maximum over $i\in[m]$, we obtain
	\begin{align*}
		\max_{i\in[m]}\Bigl(F_i(x^k+\gamma d^k)-F_i(x^k)\Bigr)
		\le& \gamma \max_{i\in[m]}
		\Bigl(g_i(x^k+d^k)-g_i(x^k)+\inner{\nabla f_i(x^k)}{d^k}\Bigr)\\
		&+\gamma\bigl(h(x^k+d^k)-h(x^k)\bigr)+\frac{L}{2}\gamma^2\|d^k\|^2 \\
		= &\gamma \max_{i\in[m]}\phi_i(x^k+d^k;x^k)
		+\gamma\bigl(h(x^k+d^k)-h(x^k)\bigr)+\frac{L}{2}\gamma^2\|d^k\|^2.
	\end{align*}
	By \ref{def:AOC}, 
	$
	\max_{i\in[m]}\phi_i(x^k+d^k;x^k)
	\le \sum_{i=1}^m\lambda_i^k\phi_i(x^k+d^k;x^k)+\epsilon_2^k.
	$
	Hence,
	\begin{align*}
		\max_{i\in[m]}\Bigl(F_i(x^k+\gamma d^k)-F_i(x^k)\Bigr)
		&\le \gamma \sum_{i=1}^m\lambda_i^k\phi_i(x^k+d^k;x^k)
		+\gamma\epsilon_2^k\\
		&\quad+\gamma\bigl(h(x^k+d^k)-h(x^k)\bigr)+\frac{L}{2}\gamma^2\|d^k\|^2 \\
		&= \gamma\Omega^k
		-\gamma\left(\frac{\tau}{2\beta_k}\|d^k\|^2+\epsilon_1^k\right)
		+\frac{L}{2}\gamma^2\|d^k\|^2.
	\end{align*}
	Therefore, whenever
	$
	\frac{L}{2}\gamma^2\|d^k\|^2
	\le
	\gamma\left(\frac{\tau}{2\beta_k}\|d^k\|^2+\epsilon_1^k\right)$,
	it follows that
	\[\max_{i\in[m]}\Bigl(F_i(x^k+\gamma d^k)-F_i(x^k)\Bigr)\le \gamma\Omega^k.\] 
	Since $\epsilon_1^k\ge 0$,
	the above inequality is guaranteed by $\gamma\le \tau/(\beta_{\max}L)$. Thus, every $\gamma\le \tau/(\beta_{\max}L)$ satisfies the line search condition, and hence the backtracking procedure is well defined.
	
	Set $a=\tau/(\beta_{\max}L)$ when $L>0$.
	If $a\geq1$, the initial trial is accepted.
	If $a<1$, either the initial trial is accepted or
	the last rejected trial is larger than $a$.
	In the latter case, the accepted step is larger than $a/2$.
	Thus $\gamma_k\geq\min\{1,\tau/(2\beta_{\max}L)\}$.
	When $L=0$, the initial trial is accepted and
	$\gamma_{\min}=1$.
\end{proof}
\subsection{Inexact projected gradient method for dual problem}

An important issue is how to obtain an inexact solution satisfying \ref{def:AOC} and the stopping criterion~\eqref{eq:xi-ubd}. A natural way is to solve the dual problem~\eqref{eq:dual-1}. However, two issues must be addressed. First, we need to analyze the differentiability of \eqref{eq:dual-1}. Second, the evaluation of its gradient requires the computation of the proximal operator, which may only be solved approximately. The model conditions stated in~\eqref{def:MOP} are sufficient for the smoothness of the dual objective. In particular, the common term $h$ cancels when comparing weights, so it need not be Lipschitz continuous.

We first establish the differentiability of the dual objective
and give an explicit representation of its gradient, as in value function differentiation
\citep{Danskin1967,TanabeFukudaYamashita2023Accelerated}.
\begin{theorem}
    Let $x\in\dom F$ and $\beta>0$. The function $\Phi_d(\cdot;x)$ is continuously differentiable along $\Delta_m$, with gradient represented by
    \begin{equation}
        \nabla\Phi_d(\lambda;x)=\phi(u;x),
        \qquad u=\prox_{\beta(h+\sum_{i=1}^{m}\lambda_i g_i)}
        \left(x-\beta\sum_{i=1}^{m}\lambda_i\nabla f_i(x)\right).
    \end{equation}
\end{theorem}
\begin{proof}
    For each $\lambda\in\Delta_m$, the proper, lower semicontinuous and strongly convex function
    $\Phi(\cdot,\lambda;x)$ has a unique minimizer $u(\lambda)\in\dom h$.
    Let $u=u(\lambda)$ and $v=u(\lambda')$. Adding the two strong convexity inequalities at these minimizers gives
    \[
        \frac1\beta\norm{u-v}^2
        \leq\inner{\phi(v;x)-\phi(u;x)}{\lambda-\lambda'}
        \leq\kappa_x\norm{u-v}\norm{\lambda-\lambda'},
    \]
    where $\kappa_x=(\sum_i(L_i+\norm{\nabla f_i(x)})^2)^{1/2}$.
    Thus $u(\lambda)$ is continuous. Since $\Phi(u,\cdot;x)$ is affine with gradient $\phi(u;x)$,
    Danskin's theorem~\citep{Danskin1967} gives the gradient representation
    $\nabla\Phi_d(\lambda;x)=\phi(u(\lambda);x)$ along $\Delta_m$.
    Indeed, the optimality of $u=u(\lambda)$ and $v=u(\lambda')$ gives
    \[
        \inner{\phi(v;x)}{\lambda'-\lambda}
        \leq\Phi_d(\lambda';x)-\Phi_d(\lambda;x)
        \leq\inner{\phi(u;x)}{\lambda'-\lambda}.
    \]
    The difference between the two bounds is at most
    $\beta\kappa_x^2\norm{\lambda'-\lambda}^2$.
    This proves differentiability along the simplex and continuity of the stated gradient.
\end{proof}
The preceding theorem shows that the projected gradient direction is well defined. We next give the Lipschitz bounds used in its analysis.
\begin{lemma}\label{lem:dual-smooth}
	Let $x\in\operatorname{dom}F$ and $\beta>0$. For each $i\in[m]$, let $L_i$ be a Lipschitz constant of $g_i$, and define
	\[
	\kappa_x
	:=
	(
	\sum_{i=1}^{m}
	\bigl(L_i+\|\nabla f_i(x)\|\bigr)^2
	)^{1/2}.
	\]
	For $\lambda\in\Delta_m$, let
	$
	u(\lambda)
	=
	\operatorname{prox}_{\beta(h+\sum_{i=1}^{m}\lambda_i g_i)}
	\left(
	x-\beta\sum_{i=1}^{m}\lambda_i\nabla f_i(x)
	\right)$.
	Then, for any $\lambda,\lambda'\in\Delta_m$, we have
	$\|u(\lambda)-u(\lambda')\|
	\leq
	\beta\kappa_x\|\lambda-\lambda'\|$.
	Consequently,
	\[
	\|\nabla\Phi_d(\lambda;x)-\nabla\Phi_d(\lambda';x)\|
	\leq
	\beta\kappa_x^2\|\lambda-\lambda'\|.
	\]
	In particular, $\Phi_d(\cdot;x)$ has a Lipschitz continuous gradient on $\Delta_m$ with Lipschitz constant $L_\lambda^x:=\beta\kappa_x^2$.
	
\end{lemma}

\begin{proof}
	The proof of the preceding theorem gives
	$\|u(\lambda)-u(\lambda')\|\leq\beta\kappa_x\|\lambda-\lambda'\|$.
	Since $\|\phi(u;x)-\phi(v;x)\|\leq\kappa_x\|u-v\|$, the gradient
	formula in the same theorem gives the second bound.
\end{proof}
The above lemma shows how the gradient of \eqref{eq:dual-1} is computed. In practice, however, evaluating
$\nabla\Phi_d(\lambda;x)=\phi(u(\lambda);x)$
requires the exact solution $u(\lambda)$.
Such an exact solution is generally not available, even when the proximal mappings of the individual functions $g_i$ can be computed efficiently. Consequently, the dual gradient itself can only be evaluated inexactly.

We therefore replace $u(\lambda)$ by an approximate solution $\tilde u$ satisfying an inexact optimality condition. The next result quantifies how the resulting error in the primal subproblem affects the dual gradient.
\begin{lemma}\label{lem:isubdiff}
	Suppose that $\epsilon\geq 0$, $x,\tilde{u}\in\dom F$, $\tilde{w}\in\R^n $, $\tilde{\lambda}\in\Delta_m$, and satisfies
	\begin{equation*}
		x-\beta\sum_{i=1}^{m}\tilde{\lambda}_i\nabla f_i(x)+\beta \tilde{w}-\tilde{u}\in\beta\partial_{\epsilon}\Bigl(h+\sum_{i=1}^{m}\tilde{\lambda}_i g_i\Bigr)(\tilde{u}).
	\end{equation*}
	Then for all $\lambda\in\Delta_m$, we have
	\begin{equation}\label{eq:lam-epi-diff}
			\Phi_{d}(\lambda;x)-\Phi_{d}(\tilde{\lambda};x)\leq\Phi(\tilde{u},\lambda;x)-\Phi_{d}(\tilde{\lambda};x)
			\leq\inner{\phi(\tilde{u};x)}{\lambda-\tilde{\lambda}}+\frac{\beta}{2}\norm{\tilde{w}}^2+\epsilon.
	\end{equation}
	It means that $\phi(\tilde{u};x)\in\partial^+_{\tilde{\epsilon}}\Phi_{d}(\tilde{\lambda};x)$, where $\tilde{\epsilon}=\frac{\beta}{2}\norm{\tilde{w}}^2+\epsilon$.
\end{lemma}
\begin{proof}
	By the assumption, we have
	\[
	x-\beta\sum_{i=1}^{m}\tilde{\lambda}_i\nabla f_i(x)+\beta \tilde{w}-\tilde{u}
	\in
	\beta\partial_\epsilon\Big(h+\sum_{i=1}^{m}\tilde{\lambda}_i g_i\Big)(\tilde u).
	\]
	From Lemma \ref{lem:pgap}, for all $u\in\dom h$,
	$\Phi(\tilde u,\tilde\lambda;x)\le \Phi(u,\tilde\lambda;x)-\frac{1}{2\beta}\norm{\tilde{u}-u}^2+\inner{ \tilde{w}}{\tilde{u}-u}+\epsilon$.
	Let $v=\arg\min_{u\in\R^n}\Phi(u,\tilde{\lambda};x)$. Hence, we have
	\begin{align}\label{eqn:phid-epiw}
		\Phi(\tilde{u},\tilde{\lambda};x)\leq& \min_{u\in\R^n}\Phi(u,\tilde{\lambda};x)-\frac{1}{2\beta}\norm{\tilde{u}-v}^2+\inner{\tilde{w}}{\tilde{u}-v}+\epsilon\nonumber\\
		=&\Phi_{d}(\tilde{\lambda};x)-\frac{1}{2\beta}\norm{\tilde{u}-v}^2+\inner{\tilde{w}}{\tilde{u}-v}+\epsilon\nonumber\\
		\leq&\Phi_{d}(\tilde{\lambda};x)+\frac{\beta}{2}\norm{\tilde{w}}^2+\epsilon.
	\end{align}
	From the definition of $\Phi(u,\lambda;x)$, we have
	\begin{equation*}
		\Phi(\tilde{u},\lambda;x)-\Phi(\tilde{u},\tilde{\lambda};x)
		=\sum_{i=1}^{m}(\lambda_i-\tilde{\lambda}_i)\phi_i(\tilde{u};x)=\inner{\phi(\tilde{u};x)}{\lambda-\tilde{\lambda}}.
	\end{equation*}
	Hence,
	$\Phi_{d}(\lambda;x)=\min_{u\in\R^n}\Phi(u,\lambda;x)\leq\Phi(\tilde{u},\lambda;x)\leq\Phi(\tilde{u},\tilde{\lambda};x)+\inner{\phi(\tilde{u};x)}{\lambda-\tilde{\lambda}}$. Combining with \eqref{eqn:phid-epiw}, we get the conclusion.
\end{proof}
Lemma \ref{lem:isubdiff} alone is not sufficient to design a inexact projected gradient method. We need to establish an inexact counterpart of the standard descent inequality.
\begin{lemma}\label{lem:dual-idescent}
	Let $x\in\dom F$ and $\beta>0$.
	Let $\tilde u\in\dom F$, $\tilde w\in\R^n$,
	$\tilde\lambda\in\Delta_m$, and $\epsilon_1\geq 0$ satisfy
	\[
	x-\beta\sum_{i=1}^{m}\tilde\lambda_i\nabla f_i(x)
	+\beta\tilde w-\tilde u
	\in
	\beta\partial_{\epsilon_1}
	\Bigl(h+\sum_{i=1}^{m}\tilde\lambda_i g_i\Bigr)(\tilde u).
	\]
	Set
	$\tilde\epsilon_1
	:=\epsilon_1+\frac{\beta}{2}\norm{\tilde w}^2$.
	Then there exists a constant $L_\lambda^x>0$, such that
	\begin{equation}\label{eq:dual-inexact-descent}
		\Phi_d(\tilde\lambda;x)-\Phi_d(\lambda;x)
		\leq
		\inner{\phi(\tilde u;x)}{\tilde\lambda-\lambda}
		+L_\lambda^x\norm{\lambda-\tilde\lambda}^2
		+\tilde\epsilon_1,
		\quad \forall\lambda\in\Delta_m.
	\end{equation}
\end{lemma}
\begin{proof}
	We first give a quadratic lower bound for $\Phi_{d}(\cdot)$ at $\tilde u$.
	By \eqref{eqn:phid-epiw} in Lemma \ref{lem:isubdiff},
	\[
	0\leq
	\Phi(\tilde u,\tilde\lambda;x)-\Phi_d(\tilde\lambda;x)
	\leq \tilde\epsilon_1.
	\]
	For any $v\in\dom F$, the $1/\beta$-strong convexity of
	$\Phi(\cdot,\tilde\lambda;x)$, applied at the midpoint
	of $v$ and $\tilde u$, gives
	\[
	\Phi_d(\tilde\lambda;x)
	\leq
	\Phi(\frac{v+\tilde u}{2},\tilde\lambda;x)
	\leq
	\frac12\Phi(v,\tilde\lambda;x)
	+\frac12\Phi(\tilde u,\tilde\lambda;x)
	-\frac{1}{8\beta}\norm{v-\tilde u}^2.
	\]
	Rearranging and using the subproblem error bound above,
	we obtain
	\begin{align}\label{eq:dual-inexact-growth}
		\Phi(v,\tilde\lambda;x)
		&\geq
		2\Phi_d(\tilde\lambda;x)
		-\Phi(\tilde u,\tilde\lambda;x)
		+\frac{1}{4\beta}\norm{v-\tilde u}^2\nonumber\\
		&\geq
		\Phi_d(\tilde\lambda;x)
		+\frac{1}{4\beta}\norm{v-\tilde u}^2
		-\tilde\epsilon_1.
	\end{align}
	
	For $u,v\in\dom h$ and $\lambda,\lambda'\in\Delta_m$, the definition of $\phi_i$ and the Lipschitz continuity of $g_i$ give
    \[
        \left|\inner{\phi(u;x)-\phi(v;x)}{\lambda-\lambda'}\right|
        \leq\sum_{i=1}^{m}|\lambda_i-\lambda_i'|(L_i+\norm{\nabla f_i(x)})\norm{u-v}
        \leq\kappa_x\norm{u-v}\norm{\lambda-\lambda'},
    \]
    where $\kappa_x$ is defined in Lemma~\ref{lem:dual-smooth}.
    The common term $h$ does not enter this estimate.

	Now fix $\lambda\in\Delta_m$ and let $v=u(\lambda)$
	be the exact minimizer of $\Phi(\cdot,\lambda;x)$.
	Since $\Phi(v,\cdot;x)$ is affine, we can combine
	\eqref{eq:dual-inexact-growth} to obtain
	\begin{align*}
		\Phi_d(\lambda;x)
		&=
		\Phi(v,\tilde\lambda;x)
		+\inner{\phi(v;x)}{\lambda-\tilde\lambda}\\
		&\geq
		\Phi_d(\tilde\lambda;x)
		+\inner{\phi(\tilde u;x)}{\lambda-\tilde\lambda}
		+\frac{\norm{v-\tilde u}^2}{4\beta}
		-\kappa_x \norm{v-\tilde u}\cdot\norm{\lambda-\tilde\lambda}
		-\tilde\epsilon_1.
	\end{align*}
	The positive quadratic term absorbs the mixed term:
	\[\frac{\norm{v-\tilde u}^2}{4\beta}-\kappa_x \norm{v-\tilde u}\norm{\lambda-\tilde\lambda}
	=\frac{1}{4\beta}
	\bigl(\norm{v-\tilde u}-2\beta\kappa_x\norm{\lambda-\tilde\lambda}\bigr)^2
	-\beta\kappa_x^2\norm{\lambda-\tilde\lambda}^2.\]
	Hence the left side above is bounded below by $-\beta\kappa_x^2\norm{\lambda-\tilde\lambda}^2$.
	Substituting this bound and rearranging proves
	\eqref{eq:dual-inexact-descent} with any
	$L_\lambda^x>0$ satisfying
	$L_\lambda^x\geq\beta\kappa_x^2$.
\end{proof}
The two inexact estimates identify the error that
must be controlled when approximate solutions are
used in the iteration. We now formulate the method in
terms of a solver whose output satisfies a prescribed
accuracy requirement.
For fixed $x\in\operatorname{dom}F$ and $\beta>0$, let $(u,w,\epsilon_1)=\mathcal{A}(x,\beta,\lambda,\alpha)$, where $\mathcal{A}$ denote an inexact proximal solver which, for
$\lambda\in\Delta_m$ and $\alpha>0$, returns a tuple satisfying $u,u-\beta w\in\dom h$ and 
\begin{equation}\label{eq:inner-A}
	\left\{
	\begin{aligned}
		&x-\beta\sum_{i=1}^{m}\lambda_i\nabla f_i(x)
		+\beta w-u
		\in
		\beta\partial_{\epsilon_1}
		\Bigl(h+\sum_{i=1}^{m}\lambda_i g_i\Bigr)(u),\\
		&\epsilon_1+\frac{\beta}{2}\|w\|^2\leq\alpha.
	\end{aligned}
	\right.
\end{equation}
The tolerance controls the error of the estimates. The stopping test must additionally account for the residual in the prescribed stopping criterion.
Define
\begin{align*}
	\mathcal H(u,w,\lambda)
	={}&
	\frac{\tau}{2\beta}
	\left(
	\|u-\beta w-x\|^2-\|u-x\|^2
	\right)
	+\frac{\beta}{2}\|w\|^2\\
	&+
	\sum_{i=1}^{m}\lambda_i
	\left(
	g_i(u-\beta w)-g_i(u)
	-\beta\langle\nabla f_i(x),w\rangle
	\right)\\
	&+h(u-\beta w)-h(u).
\end{align*}
At each inner iteration, define
\begin{align*}
	\tilde\epsilon_1^t
	=\epsilon_1^t+\frac{\beta}{2}\|w^t\|^2,\quad 
	\epsilon_2^t=\max_{i\in[m]}\phi_i(u^t-\beta w^t;x)
	-\sum_{i=1}^m\lambda_i^t\phi_i(u^t-\beta w^t;x),
\end{align*}
and $\xi^t=\mathcal H(u^t,w^t,\lambda^t)+2\epsilon_1^t+\epsilon_2^t$.

\paragraph{\textbf{Line search}}
Initialize $\gamma^+=\min\{c_\tau\gamma_t,\gamma_{\max}\}$ and compute
$
\lambda^+
=\Pi_{\Delta_m}\bigl(\lambda^t+\gamma^+\phi(u^t;x)\bigr)$.
On the first trial, if $\norm{\lambda^+-\lambda^t}^2\leq\gamma^+\alpha_t$, return \texttt{refine} to Algorithm~\ref{algo:projected}. Otherwise, compute
\[
\alpha^+
=\frac{1-\tau_\lambda}{8c_\tau\gamma^+}
\|\lambda^+-\lambda^t\|^2,
\qquad
(u^+,w^+,\epsilon_1^+)
=\mathcal A(x,\beta,\lambda^+,\alpha^+).
\]
Accept the trial if
\begin{equation}\label{eq:inner-backtracking}
		\Phi(u^t,\lambda^t;x)-\Phi(u^+,\lambda^+;x)
		\leq
		\langle\phi(u^t;x),\lambda^t-\lambda^+\rangle
		+\frac{\tau_\lambda}{2\gamma^+}
		\|\lambda^+-\lambda^t\|^2
		+2\tilde\epsilon_1^t.
\end{equation}
Otherwise, set $\gamma^+\leftarrow\gamma^+/2$, recompute $\lambda^+$, and repeat the trial starting with the computation of $\alpha^+$.
\begin{algorithm}[htbp]
	\caption{Inexact projected gradient method}
	\label{algo:projected}
	\small
	\begin{algorithmic}
		\REQUIRE
		$x\in\dom F$, $\lambda^0\in\Delta_m$,
		$\beta,\alpha_0,\eta>0$, $0<\gamma_0\leq\gamma_{\max}<\infty$,
		$\tau,\tau_\lambda\in(0,1)$,
		$\tau_0\in[0,1)$, $c_\tau>1$,
		and an inexact solver $\mathcal A$.
		
		\STATE Set $t=0$ and compute
		$(u^0,w^0,\epsilon_1^0)
		=\mathcal A(x,\beta,\lambda^0,\alpha_0)$.
		
		\WHILE{true}
		\STATE Compute
		$(\tilde\epsilon_1^t,\epsilon_2^t,\xi^t)$.
		
		\IF{$\xi^t\leq
			\frac{(1-\tau)\tau_0}{2\beta}\|u^t-x\|^2+\eta$}
		\RETURN
		$(u^t,w^t,\lambda^t,
		(\epsilon_1^t,\epsilon_2^t))$.
		\ENDIF
		
		\STATE Apply the line search with
		$(x,\beta,\lambda^t,u^t,\tilde\epsilon_1^t,\alpha_t,\gamma_t)$.
		
		\IF{the line search returns \texttt{refine}}
		\STATE Set $\alpha_t\leftarrow\alpha_t/2$.
		\STATE Compute
		$(u^t,w^t,\epsilon_1^t)
		=\mathcal A(x,\beta,\lambda^t,\alpha_t)$.
		\ELSE
		\STATE Set
		$
		(u^{t+1},w^{t+1},\lambda^{t+1},
		\epsilon_1^{t+1},\gamma_{t+1},\alpha_{t+1})
		=
		(u^+,w^+,\lambda^+,
		\epsilon_1^+,\gamma^+,\alpha^+)$.
		
		\STATE Set $t\leftarrow t+1$.
		\ENDIF
		\ENDWHILE
	\end{algorithmic}
\end{algorithm}

\begin{remark}\label{rem:refine}
	The \texttt{refine} in the line search is intended to ensure the validity of the line search. Since our gradient computations for the dual subproblem are inexact, the line search may fail when the norm of the difference between $\lambda^t$ and its projection is small. However, this does not necessarily mean that $\lambda^t$ is close to optimal; it may also be caused by computational errors. The \texttt{refine} checks the stopping criterion, if it is satisfied, Algorithm \ref{algo:projected} returns the current solution; otherwise, it indicates that the gradient error at the current point is too large, and a more accurate gradient needs to be recomputed.
\end{remark}

We next analyze the inner iteration. We first show that the line search is well defined and that the accepted step sizes are bounded away from zero. We then combine the projection inequality with the inexact dual estimates to derive a bound on the progress of the dual objective.
\begin{lemma}\label{lem:inner-linesearch}
	Suppose that $\mathcal A$ is given by Algorithm~\ref{algo:admm} and $\lambda^{t+1}\neq\lambda^t$. The step size of line search is well defined, and we have
	\begin{equation}
		\gamma_{\min}:=\min\left\{\gamma_0,\frac{\tau_\lambda}{4L^x_\lambda}\right\}
		\leq\gamma_{t+1}\leq\gamma_{\max}.
	\end{equation}
\end{lemma}
\begin{proof}
	First, we prove that the line search terminates after finitely many steps. For $t\geq 0$, from Lemma \ref{lem:isubdiff}, We know that $\Phi(u^{t},\lambda^{t};x)-\Phi_d(\lambda^{t};x)\leq\tilde{\epsilon}^{t}_1
	$. Hence,
	\begin{equation}
		\label{eq:phi-inexact-lbd}
		\Phi_d(\lambda^t;x)-\Phi_d(\lambda^{+};x)
		\geq \Phi(u^{t},\lambda^{t};x)-\Phi(u^{+},\lambda^{+};x)-\tilde{\epsilon}^{t}_1.
	\end{equation}
	Since $\Phi_d(\lambda;x)$ is $L_\lambda$-smooth, from Lemma \ref{lem:dual-idescent}, we have
	\begin{equation*}
		\Phi_d(\lambda^t;x)-\Phi_d(\lambda^{+};x)
		\leq \inner{\phi(u^t;x)}{\lambda^{t}-\lambda^{+}}+L_\lambda^x\norm{\lambda^{+}-\lambda^{t}}^2+\tilde{\epsilon}_1^t.
	\end{equation*}
	Hence for all $\gamma\in (0,\tau_\lambda/(2L_\lambda^x))$,
	\begin{equation*}
		\Phi_d(\lambda^t;x)-\Phi_d(\lambda^{+};x)
		\leq\inner{\phi(u^t;x)}{\lambda^{t}-\lambda^{+}}+\frac{\tau_\lambda}{2\gamma}\norm{\lambda^{+}-\lambda^{t}}^2+\tilde{\epsilon}_1^t.
	\end{equation*}
	It follows from \eqref{eq:phi-inexact-lbd} that for all $\gamma\in (0,\tau_\lambda/(2L_\lambda^x))$,
	\begin{align*}
		\Phi(u^{t},\lambda^{t};x)-\Phi(u^{+},\lambda^{+};x)
		\leq \inner{\phi(u^t;x)}{\lambda^{t}-\lambda^{+}}+\frac{\tau_{\lambda}}{2\gamma}\norm{\lambda^{+}-\lambda^{t}}^2+2\tilde{\epsilon}_1^t.
	\end{align*}
	The initial trial satisfies
	$\gamma_t\leq\min\{c_\tau\gamma_t,\gamma_{\max}\}\leq\gamma_{\max}$.
	If no halving is needed, the accepted step is at least
	$\gamma_t$. Otherwise, the preceding trial failed, so the
	accepted step is at least $\tau_\lambda/(4L_\lambda^x)$.
	Induction from $\gamma_0$ proves the stated bounds.
	The initialization also preserves $\gamma_{t+1}\leq c_\tau\gamma_t$.
	
	Next, we prove that \texttt{refine} terminates after finitely many steps. Let $u^\star$ solve~\eqref{eqn:primal} at $\lambda^t$, set	$\widehat\eta:=\eta+(1-\tau)\tau_0\norm{u^\star-x}^2/(2\beta)$. At the first trial after a failed stopping test, write
	$s=\lambda^+-\lambda^t$, $e=\tilde\epsilon_1^t$ and $v=u^t-\beta w^t$.
	By~\eqref{eq:inner-A},	$\beta\norm{w^t}^2/2\leq \tilde\epsilon_1^t\leq\alpha_t$.
	Lemma~\ref{lem:isubdiff} and strong convexity give
	$\norm{u^t-u(\lambda^t)}^2\leq2\beta \tilde\epsilon_1^t$. Let $R:=\max_{\lambda\in\Delta_m}\norm{u(\lambda)-x}
	+2\sqrt{2\beta\alpha_t}$. From Lemma~\ref{lem:dual-smooth} and compactness of $\Delta_m$, we know that $R$ is finite. Hence we have
	\[\norm{x-u^t+\beta w^t}\leq \norm{x-u(\lambda^t)}+\norm{u(\lambda^t)-u^t}+\beta\norm{w^t}\leq R.\]
	From the Lipschitz continuity of $g_i$, we have $\norm{\phi(u^t;x)}\leq \kappa_x\norm{u^t-x}\leq \kappa_x R $. The projection inequality of $\lambda$ gives
	\begin{align*}
		\Phi_p(u^t;x)-\Phi(u^t,\lambda^t;x)&=\max_{\lambda\in\Delta_m}\inner{\phi(u^t;x)}{\lambda-\lambda^t}\\
		&\leq\frac{1}{\gamma^+}\max_{\lambda\in\Delta_m} \inner{\lambda^+-\lambda^t}{\lambda-\lambda^+}+\inner{\phi(u^t;x)}{\lambda^+ -\lambda^t}\\
		&\leq(\sqrt{2}/\gamma_{\min}+\kappa_x R)\norm{\lambda^+-\lambda^t}.
	\end{align*}
	By~\eqref{eq:subgrad-bd}, the subgradient $q\in\partial h(u^t-\beta w^t)$ returned by Algorithm~\ref{algo:admm}
	satisfies $\norm{q}\leq R/\beta+\kappa_x$. Hence
	$h(u^t-\beta w^t)-h(u^t)\leq-\beta\inner{q}{w^t}
	\leq(R+\beta\kappa_x)\norm{w^t}$. Also, we have
	\[\max_{i\in[m]}\phi_i(u^t-\beta w^t;x)-\max_{i\in[m]}\phi_i(u^t;x)\leq \max_{i \in [m]}|\phi_i(u^t-\beta w^t;x)-\phi_i(u^t;x)|\leq \beta\kappa_x\norm{w^t}. \]
	We can give the upper bound of $\xi^t$,
	\begin{align*}
		\xi^t
		=&\frac{\tau}{2\beta}(\norm{u^t-\beta w^t-x}^2-\norm{u^t-x}^2)+\frac{\beta}{2}\norm{w^t}^2+h(u^t-\beta w^t)-h(u^t)+2\epsilon^t_1\\
		&+\sum_{i=1}^{m}\lambda^t_i(g_i(u^t-\beta w^t)-g_i(u^t)-\beta\inner{\nabla f_i(x)}{w^t})+\epsilon^t_2\\
		\leq& (2R+\beta\kappa_x)\norm{w^t}+2\tilde\epsilon_1^t+\Phi_p(u^t;x)-\Phi(u^t,\lambda^t;x)+\max_{i\in[m]}\phi_i(u^t-\beta w^t;x)-\max_{i\in[m]}\phi_i(u^t;x)\\
		\leq & 2(R+\beta\kappa_x)\sqrt{\frac{2\tilde\epsilon_1^t}{\beta}}+2\tilde\epsilon_1^t+(\frac{\sqrt{2}}{\gamma_{\min}}+\kappa_x R)\norm{\lambda^+-\lambda^t}.
	\end{align*}
	
	Since the stopping test failed,
	$\norm{u^\star-x}^2\leq2\norm{u^t-x}^2+2\norm{u^t-u^\star}^2$
	gives $\widehat\eta\leq2\xi^t+\norm{u^t-u^\star}^2/\beta$.
	As $\tilde\epsilon_1^t$ is bounded and $\gamma^+\leq\gamma_{\max}$, it follows that
	\[\widehat\eta\leq4(R+\beta\kappa_x)\norm{w^t}+6\tilde\epsilon_1^t+4\Phi_p(u^t;x)-4\Phi(u^t,\lambda^t;x).\]
	Then we have
	\begin{align*}
		\widehat\eta^2
		&\leq 48(R+\beta\kappa_x)^2\norm{w^t}^2+108(\tilde\epsilon_1^t)^2+48(\Phi_p(u^t;x)-\Phi(u^t,\lambda^t;x))^2\\
		&\leq48(\frac{\sqrt{2}}{\gamma_{\min}}+\kappa_x R)^2\frac{\norm{\lambda^+-\lambda^t}^2}{\gamma^+}+(\frac{96}{\beta}(R+\beta\kappa_x)^2+108\alpha^t)\tilde\epsilon_1^t.
	\end{align*}
	Since $\alpha_t\leq\bar \alpha:=\max\{\alpha_0,(1-\tau_\lambda)/(4c_\tau\gamma_{\min})\}$, we have
	\begin{equation}\label{eq:inner-refine-estimate}
		\widehat\eta^2\leq C(\frac{\norm{\lambda^+-\lambda^t}^2}{\gamma^+}+\tilde\epsilon_1^t),
	\end{equation}
	where $C=\max\{48({\sqrt{2}}/{\gamma_{\min}}+\kappa_x R)^2,({96}(R+\beta\kappa_x)^2/{\beta}+108\bar\alpha)\}$.
	At a \texttt{refine} both terms in parentheses are at most $\bar\alpha$.
	Thus $\bar\alpha\geq\widehat\eta^2/(2C)$ whenever \texttt{refine} is triggered.
	Successive halvings can meet this lower bound only finitely many times,
	since $\widehat\eta>0$.
\end{proof}

\begin{lemma}\label{lem:lam-idescent}
	Denote $\tilde{\epsilon}^t_1=\frac{\beta}{2}\norm{w^t}^2+\epsilon_1^t$. For all $\lambda\in\Delta_m$, we have
	\begin{equation}\label{eq:lam-iconv}
		\begin{aligned}
		&\Phi_d(\lambda;x)-\Phi_d(\lambda^{t+1};x)\\
		\leq&\frac{1}{2\gamma_{t+1}}\left(\norm{\lambda^{t}-\lambda}^2-\norm{\lambda^{t+1}-\lambda}^2-(1-\tau_{\lambda})\norm{\lambda^{t+1}-\lambda^t}^2\right)+3\tilde{\epsilon}^t_1+\tilde{\epsilon}_1^{t+1}.
		\end{aligned}
	\end{equation}
\end{lemma}
\begin{proof}
	For $t\geq 0$, from Lemma \ref{lem:isubdiff}, we know that $\Phi(u^{t},\lambda^{t};x)-\Phi_d(\lambda^{t};x)\leq\tilde{\epsilon}^{t}_1$.
	Hence,
	\begin{equation*}
		\Phi_d(\lambda^t;x)-\Phi_d(\lambda^{t+1};x)
		\leq \Phi(u^{t},\lambda^{t};x)-\Phi(u^{t+1},\lambda^{t+1};x)+\tilde{\epsilon}^{t+1}_1.
	\end{equation*}
	Combining with line search, it implies that
	\begin{equation}\label{eq:lam-descent}
		\Phi_d(\lambda^t;x)-\Phi_d(\lambda^{t+1};x)
		\leq \inner{\phi(u^t;x)}{\lambda^{t}-\lambda^{t+1}}+\frac{\tau_{\lambda}}{2\gamma_{t+1}}\norm{\lambda^{t+1}-\lambda^{t}}^2+2\tilde{\epsilon}_1^t+\tilde{\epsilon}_1^{t+1}.
	\end{equation}
	Since
	\[x-\beta\sum_{i=1}^{m}\lambda^t_i\nabla f_i(x)+\beta w^t-u^t\in\beta\partial_{\epsilon^t_1}\Bigl(h+\sum_{i=1}^{m}\lambda^t_i g_i\Bigr)(u^t),\]
	it follows from Lemma \ref{lem:isubdiff} that, for all $\lambda\in\Delta_m$, we have
	\begin{equation}\label{eq:lam-iconvex}
		\Phi_{d}(\lambda;x)-\Phi_{d}(\lambda^t;x)\leq\inner{\phi(u^t;x)}{\lambda-\lambda^t}+\tilde{\epsilon}^t_1.
	\end{equation}
	Combining with \eqref{eq:lam-descent} and \eqref{eq:lam-iconvex}, we have
	\begin{equation}\label{eq:lam-idescent-convex}
		\Phi_{d}(\lambda;x)-\Phi_{d}(\lambda^{t+1};x)
		\leq\inner{\phi(u^t;x)}{\lambda-\lambda^{t+1}}+\frac{\tau_{\lambda}}{2\gamma_{t+1}}\norm{\lambda^{t+1}-\lambda^t}^2+3\tilde{\epsilon}^t_1+\tilde{\epsilon}_1^{t+1}.
	\end{equation}
	By the characterization of the projection operator, for all $\lambda\in\Delta_m$, we have
	\[\norm{\lambda^{t+1}-\lambda}^2\leq\inner{\lambda^{t+1}-\lambda}{\lambda^t+\gamma_{t+1}\phi(u^t;x)-\lambda}.\]
	Therefore, $\gamma_{t+1}\inner{\phi(u^t;x)}{\lambda-\lambda^{t+1}}\leq\inner{\lambda^{t+1}-\lambda}{\lambda^t-\lambda}-\norm{\lambda^{t+1}-\lambda}^2$. Hence, from $\frac{1}{2}\norm{y-z}^2-\frac{1}{2}\norm{v-z}^2
	=\inner{v-z}{y-v}+\frac{1}{2}\norm{y-v}^2$, we have
	\begin{equation}\label{eq:lam-tri-ineq}
		\gamma_{t+1}\inner{\phi(u^t;x)}{\lambda-\lambda^{t+1}}\leq\frac{1}{2}\left(\norm{\lambda^{t}-\lambda}^2-\norm{\lambda^{t+1}-\lambda}^2-\norm{\lambda^{t+1}-\lambda^t}^2\right).
	\end{equation}
	Combining with \eqref{eq:lam-idescent-convex} and \eqref{eq:lam-tri-ineq}, we have
	\begin{align*}
		&\Phi_{d}(\lambda;x)-\Phi_{d}(\lambda^{t+1};x)\\
		\leq&\frac{1}{2\gamma_{t+1}}\left(\norm{\lambda^{t}-\lambda}^2-\norm{\lambda^{t+1}-\lambda}^2-\norm{\lambda^{t+1}-\lambda^t}^2\right)+\frac{\tau_{\lambda}}{2\gamma_{t+1}}\norm{\lambda^{t+1}-\lambda^t}^2+3\tilde{\epsilon}^t_1+\tilde{\epsilon}_1^{t+1}\\
		=&\frac{1}{2\gamma_{t+1}}\left(\norm{\lambda^{t}-\lambda}^2-\norm{\lambda^{t+1}-\lambda}^2-(1-\tau_{\lambda})\norm{\lambda^{t+1}-\lambda^t}^2\right)+3\tilde{\epsilon}^t_1+\tilde{\epsilon}_1^{t+1}.
	\end{align*}
\end{proof}
Now, we can establish the convergence rate of Algorithm \ref{algo:projected}.
\begin{theorem}\label{thm:projectd-conv}
	For all $t\ge 2$, we have
	\begin{equation}
		\min_{i=2,...,t+1}\{\Phi_d(\lambda;x)-\Phi_d(\lambda^{i};x)\}\leq\frac{2}{t\gamma^\lambda_{\min}}.
		\label{eq:final-rate}
	\end{equation}
\end{theorem}

\begin{proof}
	For $t\geq 1$, we have $\tilde{\epsilon}^t_1\leq \frac{1-\tau_{\lambda}}{8c_\tau\gamma_t}\norm{\lambda^{t}-\lambda^{t-1}}^2$ and $\tilde{\epsilon}^{t+1}_1\leq \frac{1-\tau_{\lambda}}{8c_\tau\gamma_{t+1}}\norm{\lambda^{t+1}-\lambda^t}^2$.
	Combining with Lemma \ref{lem:lam-idescent} and the definition of $\alpha_t$, we have
	\begin{align*}
		\Phi_d(\lambda;x)-\Phi_d(\lambda^{t+1};x)\leq& \frac{1}{2\gamma_{t+1}}\left(\norm{\lambda^t-\lambda}^2-\norm{\lambda^{t+1}-\lambda}^2\right)\\
		&-\frac{(1-\tau_{\lambda})(4c_\tau-1)}{8c_\tau\gamma_{t+1}}\norm{\lambda^{t+1}-\lambda^t}^2+\frac{3(1-\tau_{\lambda})}{8\gamma_{t+1}}\norm{\lambda^{t}-\lambda^{t-1}}^2.
	\end{align*}
	Since $c_\tau>1$, $(4c_\tau-1)/c_\tau>3$ holds. Then we have
	\begin{align*}
		\sum_{i=2}^{t+1}2\gamma_i\left(\Phi_d(\lambda;x)-\Phi_d(\lambda^{i};x)\right)
		\leq& \norm{\lambda^1-\lambda}^2-\norm{\lambda^{t+1}-\lambda}^2+\frac{3(1-\tau_\lambda)}4\norm{\lambda^{1}-\lambda^{0}}^2\\
		\leq&\norm{\lambda^1-\lambda}^2+\frac{3(1-\tau_\lambda)}4\norm{\lambda^{1}-\lambda^{0}}^2\leq 4,
	\end{align*}
	the last inequality follows from that $\norm{a-b}\leq \sqrt{2}$ for all $a,b\in\Delta_m$. Then we have
	\begin{equation*}
		\min_{i=2,...,t+1}\{\Phi_d(\lambda;x)-\Phi_d(\lambda^{i};x)\}\leq\frac{2}{t\gamma^\lambda_{\min}}.
	\end{equation*}
	This establishes \eqref{eq:final-rate} and completes the proof.
\end{proof}
\begin{remark}\label{rem:inner-complexity}
	With Algorithm~\ref{algo:admm} used as $\mathcal A$ and a fixed
	$\alpha_0>0$, Algorithm~\ref{algo:projected} satisfies its stopping
	criterion after $\mathcal O(1+\widehat\eta^{-1})$ calls to $\mathcal A$,
	including the initialization as well as the \texttt{refine} and
	backtracking steps. We can get this conclusion by three steps. First, the descent
	estimates \eqref{eq:lam-descent} and \eqref{eq:lam-iconv} imply that
	among the first $T$ accepted updates there is an iterate for which both
	$\xi^t$ and $\norm{u^t-u^\star}^2$ are bounded above by a multiple of
	$1/T$; they also imply that the sum of the incoming tolerances after
	update $n$ is $\mathcal O(1/n)$. Second, by
	Lemma~\ref{lem:inner-linesearch}, the $j$th \texttt{refine} requires an
	incoming tolerance of at least $2^{j-1}\widehat\eta^2/(2C)$. Combining
	this lower bound with the tail bound on the sum of incoming tolerances
	shows that, after initialization, the total weight of such
	\texttt{refine} steps is
	$\mathcal O(2^{-(j-1)/2}/\widehat\eta)$. Finally, summing this bound over
	$j$, and adding the initial halvings and the bounded number of
	backtracking trials, yields the claimed call complexity.
\end{remark}

\subsection{ADMM for the proximal operator}
The proximal mapping of $h+\sum_i\lambda_i g_i$ may not be available in closed form.
We therefore use ADMM to compute it from the individual proximal operators.
Without loss of generality, assume that $\lambda_i>0$ for all $i\in[m]$.
If some weights vanish, omit those components and use the number of retained components in place of $m$ in the ADMM updates; the original $m$ still bounds the cost. For fixed $z\in\R^n$, the weighted proximal problem is equivalent to
\begin{equation}\label{prob:admm}
    \min_{x_1,\ldots,x_m,\bar x}
    \sum_{i=1}^{m}\left(\beta\lambda_i g_i(x_i)+\frac1{2m}\norm{x_i-z}^2\right)+\beta h(\bar x),
    \qquad x_i=\bar x,\quad i\in[m].
\end{equation}
This is a two-block consensus problem \citep[Section~7.1]{Boyd2011ADMM}.
Updating the common variable first gives
\begin{equation}\label{eqn:admm-iter}
    \left\{
    \begin{aligned}
        \bar x^{l+1}
        &=\prox_{\rho\beta h/m}\left(\frac1m\sum_{i=1}^{m}(x_i^l+a_i^l)\right),\\
        x_i^{l+1}
        &=\prox_{\frac{m\rho\beta\lambda_i}{m+\rho}g_i}
        \left(\frac{m\bar x^{l+1}-ma_i^l+\rho z}{m+\rho}\right),\qquad i\in[m],\\
        a_i^{l+1}&=a_i^l+x_i^{l+1}-\bar x^{l+1},\qquad i\in[m].
    \end{aligned}
    \right.
\end{equation}
Let $\nu_i^l=\rho z+m\bar x^{l+1}-ma_i^l-(\rho+m)x_i^{l+1}$.
The component proximal conditions give
$\nu_i^l/(m\rho\beta\lambda_i)\in\partial g_i(x_i^{l+1})$.
To ensure that both the approximate point and its corrected point belong to
$\dom h$, define
\begin{equation}\label{eqn:admm-w}
    u^{l+1}=\bar x^{l+1},\qquad
    w^{l+1}=\frac1\beta\left[u^{l+1}-
    \prox_{\beta h}\left(z-\frac1{m\rho}\sum_{i=1}^{m}\nu_i^l\right)\right].
\end{equation}
By Proposition~\ref{prop:subdiff}, ${\nu_i^l}/{(m\rho\beta\lambda_i)}
\in\partial_{\epsilon^{l+1,i}}g_i(u^{l+1})$,
where
\begin{equation}\label{eqn:admm-epi}
    \epsilon^{l+1,i}=g_i(u^{l+1})-g_i(x_i^{l+1})
    -\frac1{m\rho\beta\lambda_i}\inner{\nu_i^l}{u^{l+1}-x_i^{l+1}}\geq0.
\end{equation}
The optimality condition of~\eqref{eqn:admm-w} is
\begin{equation}\label{eq:admm-h-subgrad}
    q^{l+1}:=\frac{z-u^{l+1}+\beta w^{l+1}}\beta
    -\frac1{m\rho\beta}\sum_{i=1}^{m}\nu_i^l
    \in\partial h(u^{l+1}-\beta w^{l+1}).
\end{equation}
Transporting this subgradient to $u^{l+1}$, set
\begin{equation}\label{eq:admm-total-error}
        \epsilon^{l+1}=\sum_i\lambda_i\epsilon^{l+1,i}
        +h(u^{l+1})-h(u^{l+1}-\beta w^{l+1})-\beta\inner{q^{l+1}}{w^{l+1}}.
\end{equation}
The last three terms form a nonnegative subgradient error.
Adding the corresponding inequality to the weighted inequalities yields
\[
    z+\beta w^{l+1}-u^{l+1}\in\beta\partial_{\epsilon^{l+1}}
    \Bigl(h+\sum_{i=1}^{m}\lambda_i g_i\Bigr)(u^{l+1}).
\]
With $z=x-\beta\sum_i\lambda_i\nabla f_i(x)$, this is precisely~\eqref{eq:inner-A}
once $\epsilon^{l+1}+\beta\norm{w^{l+1}}^2/2\leq\alpha$.
The residual $w^{l+1}$ need not vanish.

For later use, define
\[
    r^{l+1}=\frac1\beta\left(\bar x^{l+1}
    -\frac1m\sum_i x_i^{l+1}
    +\frac1\rho\sum_i(x_i^l-x_i^{l+1})\right).
\]
The common update gives
$q_0^{l+1}:=(\rho\beta)^{-1}\sum_i(x_i^l+a_i^l-\bar x^{l+1})
\in\partial h(u^{l+1})$.
Substitution shows that $q^{l+1}-q_0^{l+1}=w^{l+1}-r^{l+1}$.
By the subgradient inequality at $u^{l+1}$, the common-term error
in~\eqref{eq:admm-total-error} is at most
$\beta\inner{r^{l+1}-w^{l+1}}{w^{l+1}}$.
Consequently,
\begin{equation}\label{eq:admm-error-bound}
    \epsilon^{l+1}+\frac\beta2\norm{w^{l+1}}^2
    \leq\sum_i\lambda_i\epsilon^{l+1,i}
    +\frac\beta2\norm{r^{l+1}}^2
    -\frac\beta2\norm{w^{l+1}-r^{l+1}}^2.
\end{equation}

For every fixed $\rho>0$,~\eqref{eqn:admm-iter} is the standard two-block ADMM.
The component functions are finite on $\R^n$, and a point with
$x_i=\bar x\in\operatorname{ri}(\dom h)$ satisfies the relative interior constraint qualification.
The equivalent proximal problem has a proper, lower semicontinuous, strongly convex and coercive objective.
Thus it has a unique minimizer $x^\star=\prox_{\beta(h+\sum_i\lambda_i g_i)}(z)$,
and the unaugmented Lagrangian has a saddle point
\citep[Sections~5.2.3 and~5.4.2]{BoydVandenberghe2004}.
Standard ADMM convergence gives $x_i^l\to x^\star$, $\bar x^l\to x^\star$ and convergence of the scaled dual variables $a_i^l$.
Here $u^l=\bar x^l\to x^\star$ and $r^l\to0$.
Since the $g_i$ are Lipschitz continuous and their subgradients are bounded,
the component errors tend to zero. Equation~\eqref{eq:admm-error-bound}
therefore gives $\epsilon^l+\beta\norm{w^l}^2/2\to0$.
Hence every call with $\alpha>0$ terminates.


We can give the upper bound of the subgradient of $h$, which is used in Lemma \ref{lem:inner-linesearch}. The classic result of subgradient bound gives that $\norm{d}\le L$ for every $y\in\R^n$ and every $d\in\partial g(y)$ if $g(x)$ is convex and $L$-Lipschitz continuous. For each $i\in [m]$, the component proximal update in Algorithm~\ref{algo:admm} gives
$d_i=\nu_i^l/(m\rho\beta\lambda_i)\in\partial g_i(x_i^{l+1})$.
Then $\norm{d_i}\le L_i$. Substituting $z=x-\beta\sum_i\lambda_i\nabla f_i(x)$ into \eqref{eq:admm-h-subgrad} and writing $v^{l+1}=u^{l+1}-\beta w^{l+1}$, we obtain
\[
q^{l+1}=\frac{x-u^{l+1}+\beta w^{l+1}}{\beta}
-\sum_{i\in [m]}\lambda_i\bigl(\nabla f_i(x)+d_i\bigr)
\in\partial h(u^{l+1}-\beta w^{l+1}).
\]
Since $\lambda\in\Delta_m$, it follows that
\begin{equation}\label{eq:subgrad-bd}
	\norm{q^{l+1}}
	\le \frac{\norm{x-u^{l+1}+\beta w^{l+1}}}{\beta}
	+\max_{i\in[m]}\bigl(\norm{\nabla f_i(x)}+L_i\bigr).
\end{equation}
\begin{algorithm}[htbp]
    \caption{ADMM for Proximal Operator}
    \label{algo:admm}
    \begin{algorithmic}[1]
        \REQUIRE Parameters $\alpha,\beta,\rho>0$, vectors $\bar x^0,z\in\R^n$, and $\lambda\in\Delta_m\cap\R_{++}^m$.
        \STATE Set $l=0$, $\epsilon^0=+\infty$, $w^0=0$, $x_i^0=\bar x^0$ and $a_i^0=0$ for all $i\in[m]$.
        \WHILE{$\epsilon^l+\beta\norm{w^l}^2/2\geq\alpha$}
        \STATE Update $\bar x^{l+1}$, $x_i^{l+1}$ and $a_i^{l+1}$ by~\eqref{eqn:admm-iter};
        \STATE Compute $\nu_i^l=\rho z+m\bar x^{l+1}-ma_i^l-(\rho+m)x_i^{l+1}$;
        \STATE Compute $u^{l+1}$ and $w^{l+1}$ by~\eqref{eqn:admm-w};
        \STATE Compute $\epsilon^{l+1,i}$, $q^{l+1}$ and $\epsilon^{l+1}$ by~\eqref{eqn:admm-epi}, \eqref{eq:admm-h-subgrad} and~\eqref{eq:admm-total-error};
        \STATE $l\leftarrow l+1$;
        \ENDWHILE
        \RETURN $u^l,w^l,\epsilon^l$.
    \end{algorithmic}
\end{algorithm}
\begin{remark}\label{rem:proximal-complexity}
	Assume that proximal evaluations are exact and that $\rho\ge m$ is fixed
	within each call. Then, by the Douglas--Rachford residual bound of
	\citet{DavisYin2017} together with~\eqref{eq:admm-error-bound},
	Algorithm~\ref{algo:admm} reaches tolerance $\alpha$ in
	$\mathcal O(1+\alpha^{-1})$ iterations. Since each iteration uses at most
	$m$ proximal evaluations of the $g_i$ and two of $h$, a single call
	therefore uses $\mathcal O((m+2)(1+\alpha^{-1}))$ proximal evaluations.
	If the constants in these per-call bounds are uniform across calls---for
	example, if $\rho$ is common and the initial distances to the respective
	Douglas--Rachford fixed-point sets are uniformly bounded---then
	Remark~\ref{rem:inner-complexity} gives a total of
	$\mathcal O((m+2)(1+\alpha_0^{-1}+\widehat\eta^{-3}))$ evaluations.
\end{remark}

\section{Convergence analysis}
\label{sec:conv}
In this section, we obtain the convergence analysis of Algorithm \ref{algo:proximal} under convexity. Firstly, we give an important inequality.
\begin{lemma}\label{lem:bound}
	For all $z\in\dom h$ and $k\ge 0$, we have
	\begin{equation}\label{eq:descent-fun}
		\begin{aligned}
			&\min_{i\in [m]}\left(F_i(x^k)-F_i(z)\right)+\frac{1}{\gamma_k}\max_{i\in [m]}\left(F_i(x^{k+1})-F_i(x^k)\right)\\
			\leq&\frac{1}{2\beta_k\gamma_k}\norm{x^k-z}^2-\frac{1}{2\beta_k\gamma_k}\norm{x^{k+1}-z}^2+\eta_k.
		\end{aligned}
	\end{equation}
\end{lemma}
\begin{proof}
	Set $v^k=u^k-\beta_k w^k$. By the line search and convexity of $f_i$, followed by Lemma~\ref{lem:pgap}, we obtain
	\begin{align}\label{eqn:bound-ls}
		&\min_{i\in[m]}\bigl(F_i(x^k)-F_i(z)\bigr)
		+\frac{1}{\gamma_k}\max_{i\in[m]}\bigl(F_i(x^{k+1})-F_i(x^k)\bigr)\nonumber\\
		&\leq h(v^k)-h(z)+\sum_{i=1}^{m}\lambda_i^k\left(g_i(v^k)-g_i(z)
		+\inner{\nabla f_i(x^k)}{v^k-z}\right)
		+\frac{\tau}{2\beta_k}\norm{v^k-x^k}^2+\epsilon_1^k+\epsilon_2^k\nonumber\\
		&\leq\frac{1}{2\beta_k}\left(\norm{x^k-z}^2-\norm{v^k-z}^2\right)
		-\frac{1-\tau}{2\beta_k}\norm{u^k-x^k}^2+\xi^k\nonumber\\
		&\leq\frac{1}{2\beta_k}\left(\norm{x^k-z}^2-\norm{v^k-z}^2\right)+\eta_k.
	\end{align}
	Here the second inequality uses the definition of $\xi^k$ and the identity
	\[
	-\frac{1}{2\beta_k}\norm{u^k-z}^2+\inner{w^k}{u^k-z}
	=-\frac{1}{2\beta_k}\norm{v^k-z}^2+\frac{\beta_k}{2}\norm{w^k}^2,
	\]
	while the last follows from \eqref{eq:xi-ubd} and $\tau_0<1$.
	Since $x^{k+1}=(1-\gamma_k)x^k+\gamma_k v^k$ and $0<\gamma_k\leq1$, we have
	\begin{align*}
		\norm{x^{k+1}-z}^2
		&=(1-\gamma_k)\norm{x^k-z}^2+\gamma_k\norm{v^k-z}^2
		-\gamma_k(1-\gamma_k)\norm{v^k-x^k}^2\\
		&\leq(1-\gamma_k)\norm{x^k-z}^2+\gamma_k\norm{v^k-z}^2.
	\end{align*}
	Rearranging and dividing by $\gamma_k>0$ yields
	\[
	\norm{x^k-z}^2-\norm{v^k-z}^2
	\leq\frac{1}{\gamma_k}
	\left(\norm{x^k-z}^2-\norm{x^{k+1}-z}^2\right).
	\]
	Substituting this estimate into \eqref{eqn:bound-ls} proves \eqref{eq:descent-fun}.
\end{proof}
An immediate consequence of the fundamental inequality is
a componentwise bound on the increase in the objective values.
\begin{corollary}\label{coro:nonincreasing}
	$F_i(x^{k+1})\leq F_i(x^k)+ \gamma_k\eta_k$ for all $i\in[m]$.
\end{corollary}
\begin{proof}
	Setting $z=x^k$ in \eqref{eq:descent-fun}, we obtain $\max_{i\in [m]}\bigl(F_i(x^{k+1})-F_i(x^k)\bigr)\le \gamma_k\eta_k$. The conclusion directly follows.
\end{proof}
Since $\gamma_k\leq 1$, Corollary~\ref{coro:nonincreasing}
implies that the values $F_i(x^k)+\sum_{l=k}^{\infty}\eta_l$ are nonincreasing. This observation motivates the assumption introduced below. To obtain compactness of these sets, we impose the following boundedness assumption on the initial lower level set.
\begin{assumption}\label{assu:sol}
	The initial point satisfies $x^0\in\dom F$, and the lower level set $\tilde{U}$ which is defined as
	\begin{equation}
		\tilde{U}:=\{x\in\dom F:F(x)\le F(x^0)\}
	\end{equation}
	is nonempty and bounded.
\end{assumption}
The following lemma provides a compact set $U_0$ containing
all iterates and a nonempty set $U_\infty$ of reference points.
In particular, each objective component has a finite lower
bound on the region containing the iterates. For $i\in [m]$, denote these
bounds by
\[
F_i^*:=\min_{x\in U_0}F_i(x)>-\infty.
\]
The next theorem combines these bounds with the fundamental
inequality to establish convergence of the objective values
and square summability of the successive steps.
\begin{lemma}
	\label{lem:lset-nonempty}
	Suppose that Assumption \ref{assu:sol} holds and that $\{\eta_k\}$ is a summable nonnegative sequence. Define $U_k:=\{x\in\dom F:F(x)\leq F(x^k)+\sum_{l=k}^{\infty}\eta_l\}$ and $\sum_{l=0}^{\infty}\eta_l=\tilde{\eta}<\infty$. Then the intersection of $U_k$ is nonempty,
	\begin{equation}
			U_{\infty}:=\bigcap_{k=0}^{+\infty} U_k\neq \emptyset.
		\end{equation}
\end{lemma}
\begin{proof}
	We prove $U_0$ is bounded first. $\tilde{U}$ is bounded by Assumption \ref{assu:sol}. Since 
	\[U_0=\{x\in\dom F:\max_{i \in [m]}\{F_i(x)-F_i(x^0)\}\leq\tilde{\eta} \}.\]
	The boundedness of $U_0$ follows from \citep[Proposition~3.23]{RockafellarWets1998}.
	
	The nonempty property of $U_k$ is obvious. From the l.s.c property of $F(x)$, we know that $U_k$ is closed. By Corollary \ref{coro:nonincreasing}, we have $F_i(x^{k+1})\leq F_i(x^k)+ \gamma_k\eta_k$. Since $\gamma_k\leq 1$, we get
	\[F_i(x^{k+1})+\sum_{l=k+1}^{\infty}\eta_l\leq F_i(x^k)+\sum_{l=k}^{\infty}\eta_l.\]
	From the definition of $U_k$, $U_{k+1}\subseteq U_k$ holds for all $k\geq 0$. Hence $U_k$ is bounded, also is compact for all $k\geq 0$. Then we get the conclusion.
\end{proof}
\begin{theorem}
	\label{thm:conv}
	Suppose that Assumption \ref{assu:sol} holds and $\beta_k\in [\beta_{\min},\beta_{\max}]\subset (0,+\infty)$ and $\tilde{\eta}=\sum_{k=0}^{+\infty}\eta_k<\infty$. Then the sequence $\{x^k\}$ is bounded, and every accumulation
	point belongs to $U_\infty$.
	Moreover, $\{F_i(x^k)\}$ converges to a finite value
	for every $i\in[m]$, and
	$\sum_{k=0}^{\infty}\|x^{k+1}-x^k\|^2<\infty$.
\end{theorem}
\begin{proof}
	Since $x^k\in U_k\subseteq U_0$, $\{x^k\}\subset U_0$ is bounded. From the proof of Lemma \ref{lem:lset-nonempty}, we know that $\{F_i(x^k)+\sum_{l=k}^{\infty}\eta_l\}$ is monotonically nonincreasing. Since $F_i(x^k)\geq F_i^*$, $F_i(x^k)$ converges to a  finite value. From \eqref{eq:descent-fun}, we have
	\begin{equation*}
			\frac{1}{\gamma_k}\max_{i\in [m]}\left(F_i(x^{k+1})-F_i(x^k)\right)
			\leq\frac{1}{2\beta_k\gamma_k}\norm{x^k-z}^2-\frac{1}{2\beta_k\gamma_k}\norm{x^{k+1}-z}^2+\eta_k.
		\end{equation*}
	Let $z=x^k$. Define $\delta_k=\min_{i\in [m]}\left(F_i(x^{k})-F_i(x^{k+1})\right)$. Since $\gamma_k\leq 1$, we have
	\begin{equation}
			\norm{x^{k+1}-x^k}^2\leq2\beta_k\delta_k+2\beta_k\eta_k.
		\end{equation}
	We have
	\begin{align*}
		\sum_{k=0}^{+\infty}(\delta_k+\eta_k)\leq&
		\sum_{k=0}^{+\infty}\min_{i\in [m]}(F_i(x^{k})-F_i(x^{k+1}))+\tilde{\eta}
		\leq\min_{i\in [m]}\sum_{k=1}^{+\infty}(F_i(x^{k})-F_i(x^{k+1}))+\tilde{\eta}\\
		=&\min_{i\in [m]}\left(F_i(x^{0})-\lim_{k\to +\infty}F_i(x^{k+1})\right)+\tilde{\eta}\\
		\leq&\min_{i\in [m]}\left(F_i(x^{0})-F_i^*\right)+\tilde{\eta}.
	\end{align*}
	It means that 
	\[\sum_{k=0}^{+\infty}\norm{x^{k+1}-x^k}^2\leq 2\beta_{\max}\min_{i\in [m]}\left(F_i(x^{0})-F_i^*\right)+2\beta_{\max}\tilde{\eta}.\]
	The boundedness of $\{x^k\}$ implies that there is a subsequence $\{x^{k_j}\}\subseteq \{x^k\}$ such that $x^{k_j}\to x^*\in U_\infty$. For any $k\geq 0$, there is a $j_0\geq 0$ such that for any $k_j\geq k_{j_0}$, we have $x^{k_j}\in U_{k_j}\subseteq U_k$. From the closedness of $U_k$, $x^*\in U_k$ holds. From the arbitrariness of $k$, we get $x^*\in U_\infty$.
\end{proof}
Theorem~\ref{thm:conv} describes the asymptotic behavior of
the objective values.
We next establish the convergence rate of Algorithm $\ref{algo:proximal}$.
\begin{theorem}
	\label{thm:rate}
	Suppose that Assumption \ref{assu:sol} holds, $\beta_k\in [\beta_{\min},\beta_{\max}]\subset (0,+\infty)$ and $\tilde{\eta}=\sum_{k=0}^{+\infty}\eta_k$. Let $\{x^k\}$ be generated by Algorithm \ref{algo:proximal}.	Then, for each $z\in U_{\infty}$ and $k\ge 1$,
	\begin{equation}
		\label{eq:rate-phi}
		\min_{i\in [m]}\bigl(F_i(x^{k})-F_i(z)\bigr)
		\le
		\frac{\mathcal R_k(z)}{2\beta_{\min}\,\gamma_{\min}\,k},
	\end{equation}
	where
	\[\mathcal{R}_k(z)=\|x^0-z\|^2+2\beta_{\max}\min_{i\in [m]}\left(F_i(x^{0})-F_i^*\right)+2\beta_{\max}\tilde{\eta}+2\beta_{\max}\sum_{l=0}^{k-1}(l+1)\eta_l.\]
\end{theorem}

\begin{proof}
	Define $\delta_k=\min_{i\in [m]}(F_i(x^{k})-F_i(x^{k+1}))$. Multiplying \eqref{eq:descent-fun} by $\beta_k\gamma_k$ to get
	\[
	\beta_k\gamma_k\,\min_{i\in [m]}\left(F_i(x^k)-F_i(z)\right)
	\le
	\frac{1}{2}\|x^k-z\|^2-\frac{1}{2}\|x^{k+1}-z\|^2+\beta_k\delta_k+\beta_k\gamma_k\eta_k.
	\]
	Summing this inequality from $l=0$ to $k$, we obtain
	\[
	\sum_{l=0}^{k}\beta_l\gamma_l\,\min_{i\in [m]}\bigl(F_i(x^l)-F_i(z)\bigr)
	\le
	\frac{1}{2}\|x^0-z\|^2-\frac{1}{2}\|x^{k+1}-z\|^2+\sum_{l=0}^{k}(\beta_l\delta_l+\beta_l\gamma_l\eta_l).
	\]
	In Theorem \ref{thm:conv}, we have proved that $		\sum_{k=0}^{+\infty}\delta_k\leq\min_{i\in [m]}\left(F_i(x^{0})-F_i^*\right)$.	Since $\gamma_k\leq 1$ for all $k\geq 0$, we have
	\[\sum_{l=0}^{k}(\beta_l\delta_l+\beta_l\gamma_l\eta_l)\leq
	\beta_{\max}\sum_{l=0}^{k}(\delta_l+\gamma_l\eta_l)\leq \beta_{\max}\min_{i\in [m]}\left(F_i(x^{0})-F_i^*\right)+\beta_{\max}\sum_{l=0}^{k}\eta_l.\]
	Lemma \ref{lem:bound} implies
	$\min_{i\in [m]}\bigl(F_i(x^{k+1})-F_i(z)\bigr)-\min_{i\in [m]}\bigl(F_i(x^k)-F_i(z)\bigr)\leq \eta_k$.
	By Corollary~\ref{coro:nonincreasing} and $\gamma_j\leq1$,
	for $0\leq l\leq k-1$, we have
	\[
	\min_{i\in[m]}\bigl(F_i(x^k)-F_i(z)\bigr)
	\leq
	\min_{i\in[m]}\bigl(F_i(x^l)-F_i(z)\bigr)
	+\sum_{j=l}^{k-1}\eta_j.
	\]
	Multiplying by $\beta_l\gamma_l$ and summing over
	$l=0,\ldots,k-1$ yields
	\[
		\left(\sum_{l=0}^{k-1}\beta_l\gamma_l\right)
		\min_{i\in[m]}\bigl(F_i(x^k)-F_i(z)\bigr)
		\leq
		\sum_{l=0}^{k-1}\beta_l\gamma_l
		\min_{i\in[m]}\bigl(F_i(x^l)-F_i(z)\bigr)
		+\sum_{l=0}^{k-1}\beta_l\gamma_l
		\sum_{j=l}^{k-1}\eta_j.
	\]
	Changing the order of summation and using
	$\beta_l\gamma_l\leq\beta_{\max}$, we obtain
	\[
		\sum_{l=0}^{k-1}\beta_l\gamma_l
		\sum_{j=l}^{k-1}\eta_j
		=
		\sum_{j=0}^{k-1}\eta_j
		\sum_{l=0}^{j}\beta_l\gamma_l
		\leq
		\beta_{\max}\sum_{j=0}^{k-1}(j+1)\eta_j.
	\]
	Then we have
	\[
	\min_{i\in [m]}\bigl(F_i(x^{k})-F_i(z)\bigr)
	\le
	\frac{\mathcal{R}_k(z)}{2\sum_{l=0}^{k-1}\beta_l\gamma_l}.
	\]
	Finally, using $\beta_l\ge \beta_{\min}$ and $\gamma_l\ge \gamma_{\min}$ for all $l\ge 0$, we obtain
	\[
	\min_{i\in [m]}\bigl(F_i(x^{k})-F_i(z)\bigr)
	\le
	\frac{\mathcal{R}_k(z)}{2\beta_{\min}\,\gamma_{\min}\,k}.
	\]
	This completes the proof.
\end{proof}
\begin{remark}
	The error dependent terms in \eqref{eq:rate-phi}
	vanish when $\eta_k\equiv0$.
	Weighted summability,
	$\sum_{k=0}^{\infty}(k+1)\eta_k<\infty$,
	yields an $\mathcal O(1/k)$ rate.
	The same requirement appears in
	\citet[Theorem~3.1(iii)]{YangToh2026}
	for the last iterate rate under an absolute stopping
	criterion, with $\Delta^k=0$ in their notation.
	Under mere summability of $\{\eta_k\}$, the bound
	still tends to zero, but this rate is not guaranteed.
\end{remark}
Theorem \ref{thm:rate} also allows us to identify the
accumulation points. Under the assumed summability of the
errors, the upper bound in~\eqref{eq:rate-phi} tends to zero.
If an accumulation point admitted a strict improvement in
every objective component, the corresponding objective
residual would remain bounded away from zero along a
convergent subsequence. The following theorem rules out
this possibility and establishes weak Pareto optimality.
\begin{theorem}
	Suppose that Assumption \ref{assu:sol} holds and $\beta_k\in [\beta_{\min},\beta_{\max}]\subset (0,+\infty)$. Let $\{x^k\}$ be generated by Algorithm \ref{algo:proximal}. Let $x^*$ be any accumulation point of $\{x^k\}$. Then $x^*$ is a weak Pareto optimal solution of \eqref{def:MOP}.
\end{theorem}
\begin{proof}
	From Theorem \ref{thm:rate}, for any $z\in U_{\infty}$, we know that
	\begin{equation*}
		\min_{i\in [m]}\bigl(F_i(x^{k})-F_i(z)\bigr)
		\le
		\frac{\mathcal R_k(z)}{2\beta_{\min}\,\gamma_{\min}\,k}.
	\end{equation*}
	Assume that $x^*$ is not a weak Pareto optimal solution. Then there exists a point $\bar{x}$ such that $F(\bar{x})<F(x^*)$. Obviously, $\bar{x}\in U_{\infty}$ holds. Let $\epsilon^*=\min_{i\in [m]}\left(F_i(x^*)-F_i(\bar{x})\right)$ and we have $\epsilon^*>0$. Then
	\begin{align*}
		\epsilon^*-\sum_{l=k}^{\infty}\eta_l\leq&\min_{i\in [m]}\bigl(F_i(x^{k})-F_i(x^*)\bigr)+\epsilon^*\\
		=&\min_{i\in [m]}\bigl(F_i(x^{k})-F_i(x^*)\bigr)+\min_{i\in [m]}\bigl(F_i(x^{*})-F_i(\bar{x})\bigr)\\
		\leq&
		\min_{i\in [m]}\bigl(F_i(x^{k})-F_i(\bar{x})\bigr)
		\le
		\frac{\mathcal R_k(\bar{x})}{2\beta_{\min}\,\gamma_{\min}\,k}.
	\end{align*}
	Let $k\to\infty$, the left side converges to $\epsilon^*>0$, whereas the right-hand side converges to zero. This contradiction proves the claim.
\end{proof}

\section{Numerical experiments}
\label{sec:experment}
In this section, we present some numerical results using Algorithm \ref{algo:proximal}. The experiments use Python~3.12.14 on an Intel Core i7-12800HX CPU at 2.00~GHz with 32~GB of RAM. We compare IPG with
PGM~\citep{TanabeFukudaYamashita2019}, which solves its
proximal subproblems using COPT~\citep{Ge2023COPT}. PGM uses Armijo type line search with parameter of $10^{-4}$ and a backtracking factor of $0.5$. For IPG, we set
$\tau=\tau_0=0.9$ and $c_\tau=2$.
The stopping criteria are
\begin{equation}
	\label{eq:numerical-stop}
	\|d^k\|_2\leq 5\times10^{-5},\qquad
	d^k=\begin{cases}
		u^k-\beta_k w^k-x^k, & \text{IPG},\\
		y^k-x^k, & \text{PGM}.
	\end{cases}
\end{equation}
The outer iteration limit is $1000$ and $\eta_k=10^{-3}/(k+1)^2$.

\subsection{Toy problems}

For $n=30$ and $\lambda=1/300$, consider
\begin{equation}
	\label{eq:jos1-model}
	\min_{x\in\mathbb R^n}
	\left(\frac{1}{n}\|x\|_2^2+{\lambda}g_1(x),\,
	\frac{1}{n}\|x-2\boldsymbol 1\|_2^2+{\lambda}g_2(x)\right).
\end{equation}
Here $h=0$. The two regularization choices are

\[
	\ell_1\text{--}\ell_1\colon
	 g_1(x)=g_2(x)=\|x\|_1,\qquad
	\mathrm{TV}\colon\quad
	 g_1(x)=\|x\|_1,
	g_2(x)=\sum_{j=1}^{n-1}|x_{j+1}-x_j|.
\]

The shared $\ell_1$ proximal mapping is soft thresholding;
the weighted $\ell_1$--TV mapping was evaluated by
the original ADMM implementation. All methods use the same $100$
initial points sampled uniformly from $[-10,10]^{30}$. We set $\beta_k=2$.

\begin{figure}[!htbp]
	\centering
	\includegraphics[width=0.80\linewidth]{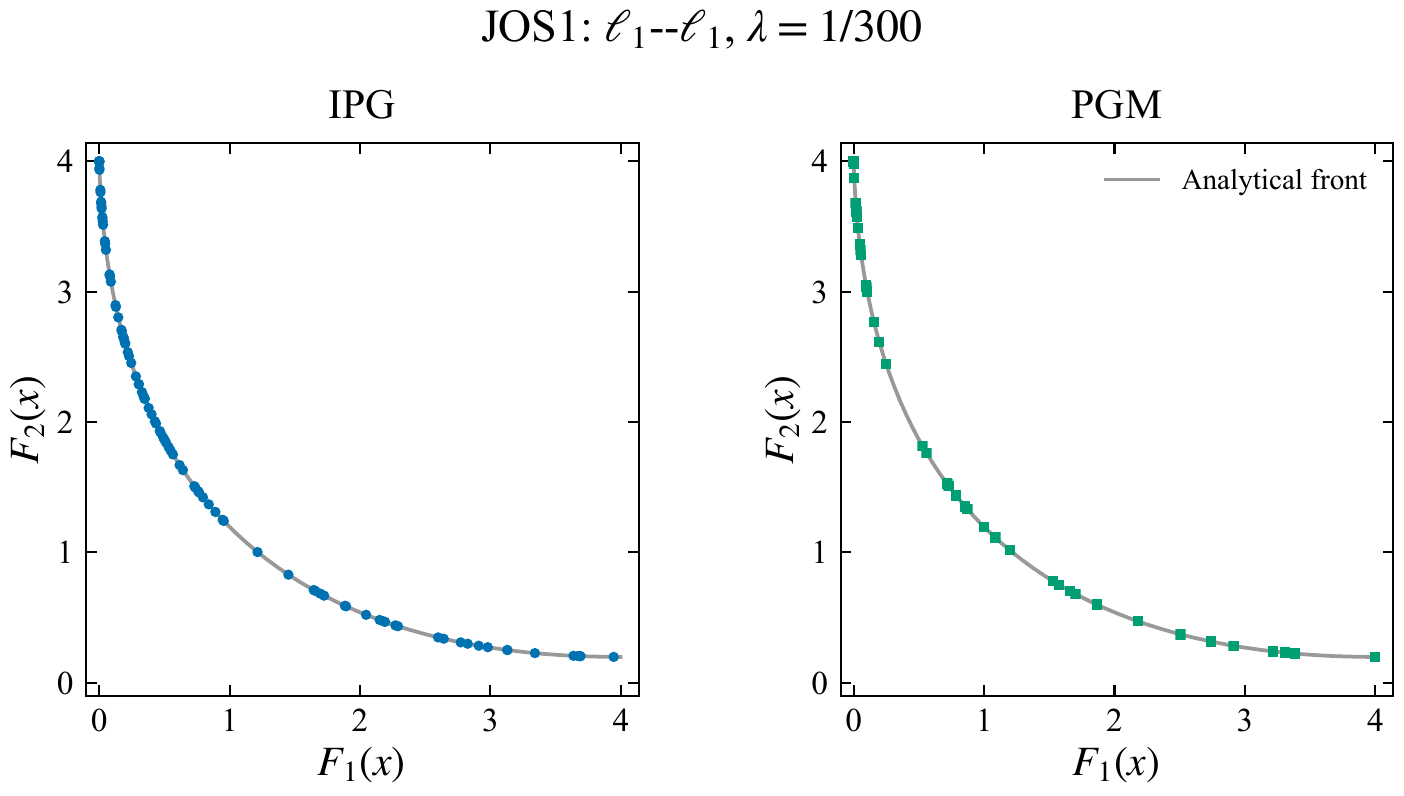}
	\caption{Pareto fronts for JOS1-$\ell_1$.}
	\label{fig:jos1-L1-fronts}
\end{figure}

\begin{figure}[!htbp]
	\centering
	\includegraphics[width=0.80\linewidth]{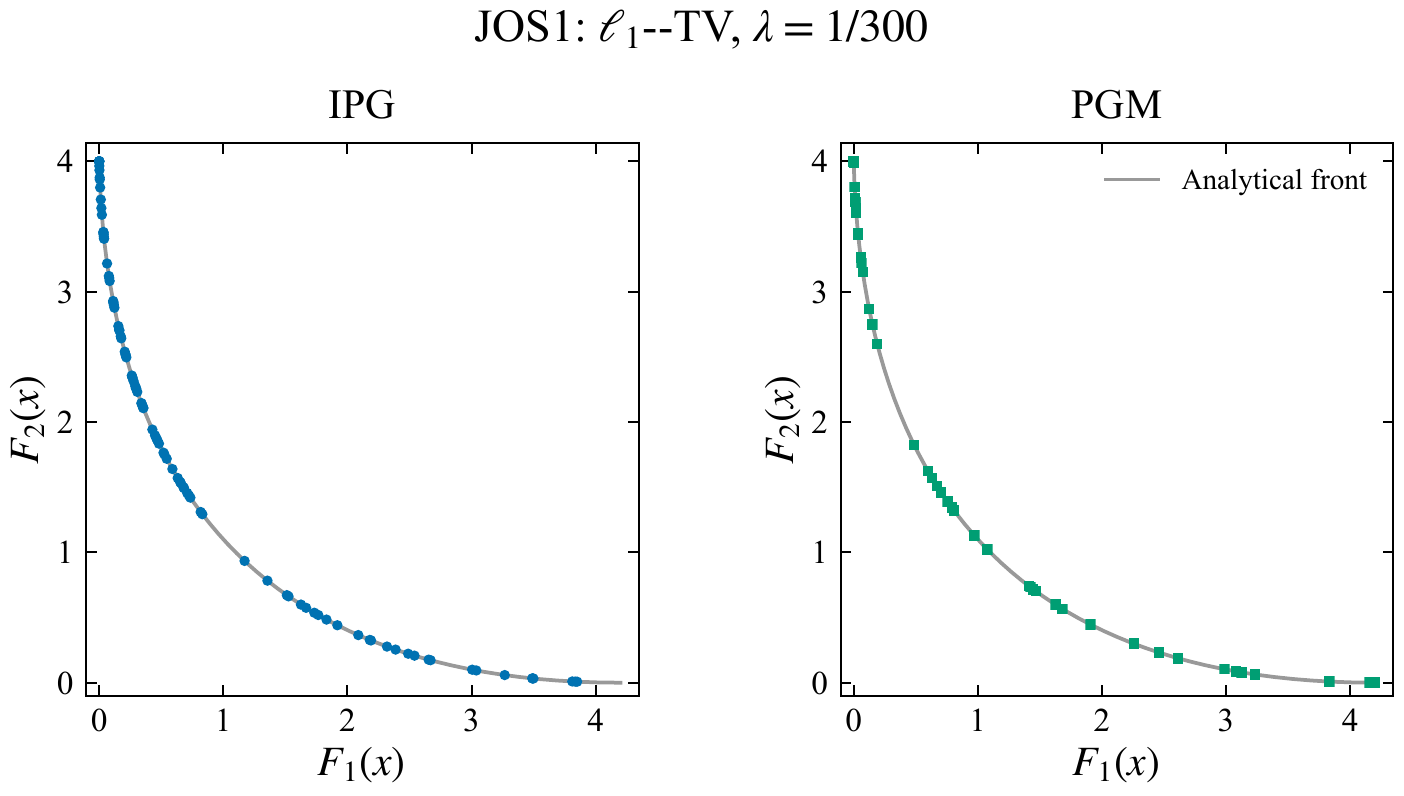}
	\caption{Pareto fronts for JOS1-$\mathrm{TV}$.}
	\label{fig:jos1-TV-fronts}
\end{figure}

\begin{table}[!htbp]
	\centering\small
	\setlength{\tabcolsep}{5pt}
	\caption{JOS1 results ($n=30$).}
	\label{tab:jos1-performance}
	\begin{tabular}{@{}llrrrr@{}}
		\toprule
		Regularizer & Method & Time (s) & Outer & Inner & ADMM \\
		\midrule
		$\ell_1$--$\ell_1$ & IPG & \textbf{0.0408} & 75.87 & 15.16 & -- \\
		& PGM & 0.5056 & \textbf{61.10} & 622.60 & -- \\
		\addlinespace
		$\ell_1$--TV & IPG & \textbf{0.4311} & 50.65 & 101.95 & 3296.21 \\
		& PGM & 0.4451 & \textbf{46.13} & 632.72 & -- \\
		\bottomrule
	\end{tabular}
\end{table}

Table~\ref{tab:jos1-performance}
shows that IPG is fastest for the shared $\ell_1$ regularizer despite using more outer updates. With TV regularization, the additional proximal computations increase its cost, and PGM is fastest. Figure~\ref{fig:jos1-inner} compares the inner errors for
$\ell_1$ and TV. IPG reports $\xi_k$, whereas PGM reports its primal--dual gap. The plots illustrate how IPG can accept inexact inner solutions while the outer iterates approach
the Pareto front.

\par\addvspace{\intextsep}
\noindent\begin{minipage}{\linewidth}
	\centering
	\captionsetup{type=figure}
	\includegraphics[width=0.47\linewidth]{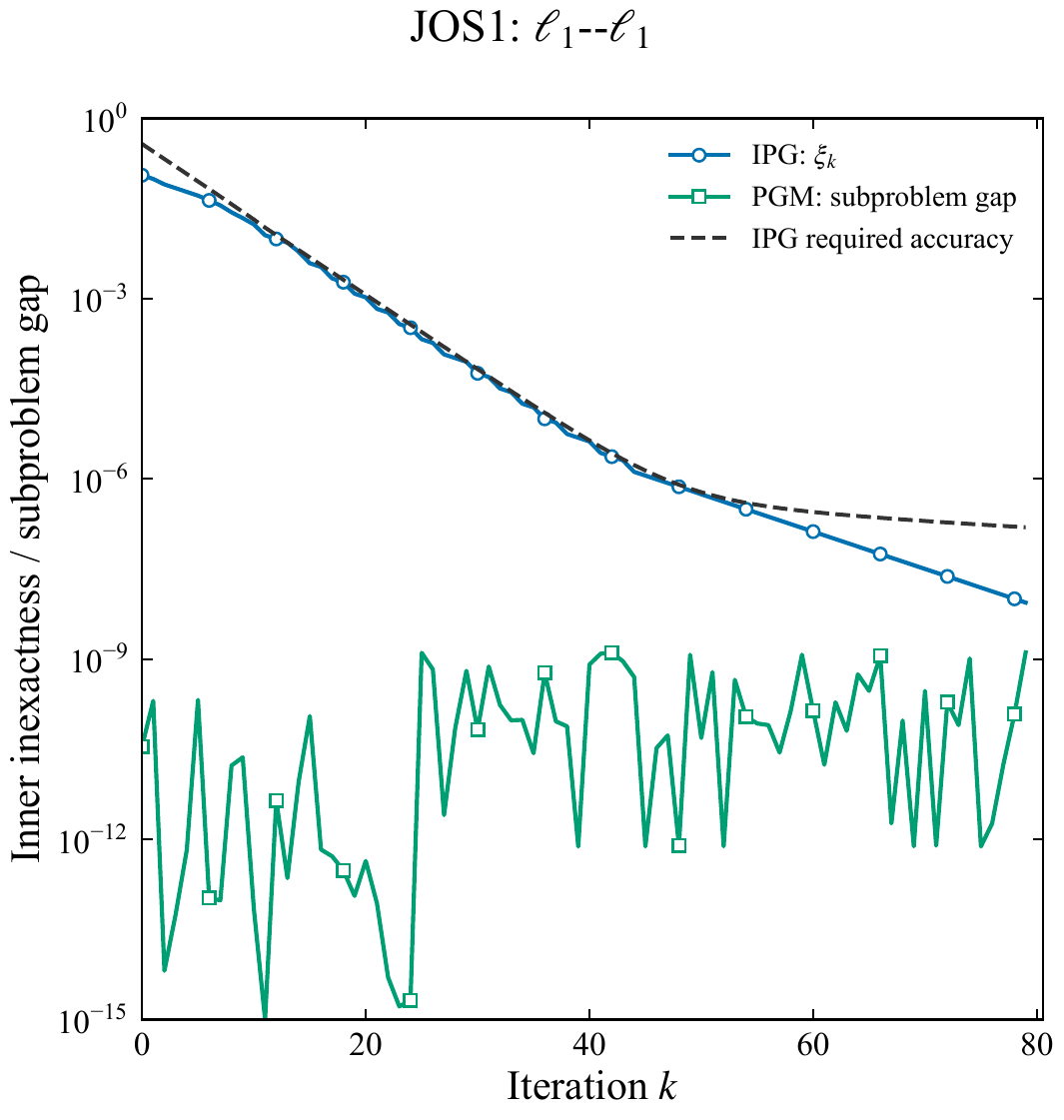}\hfill
	\includegraphics[width=0.47\linewidth]{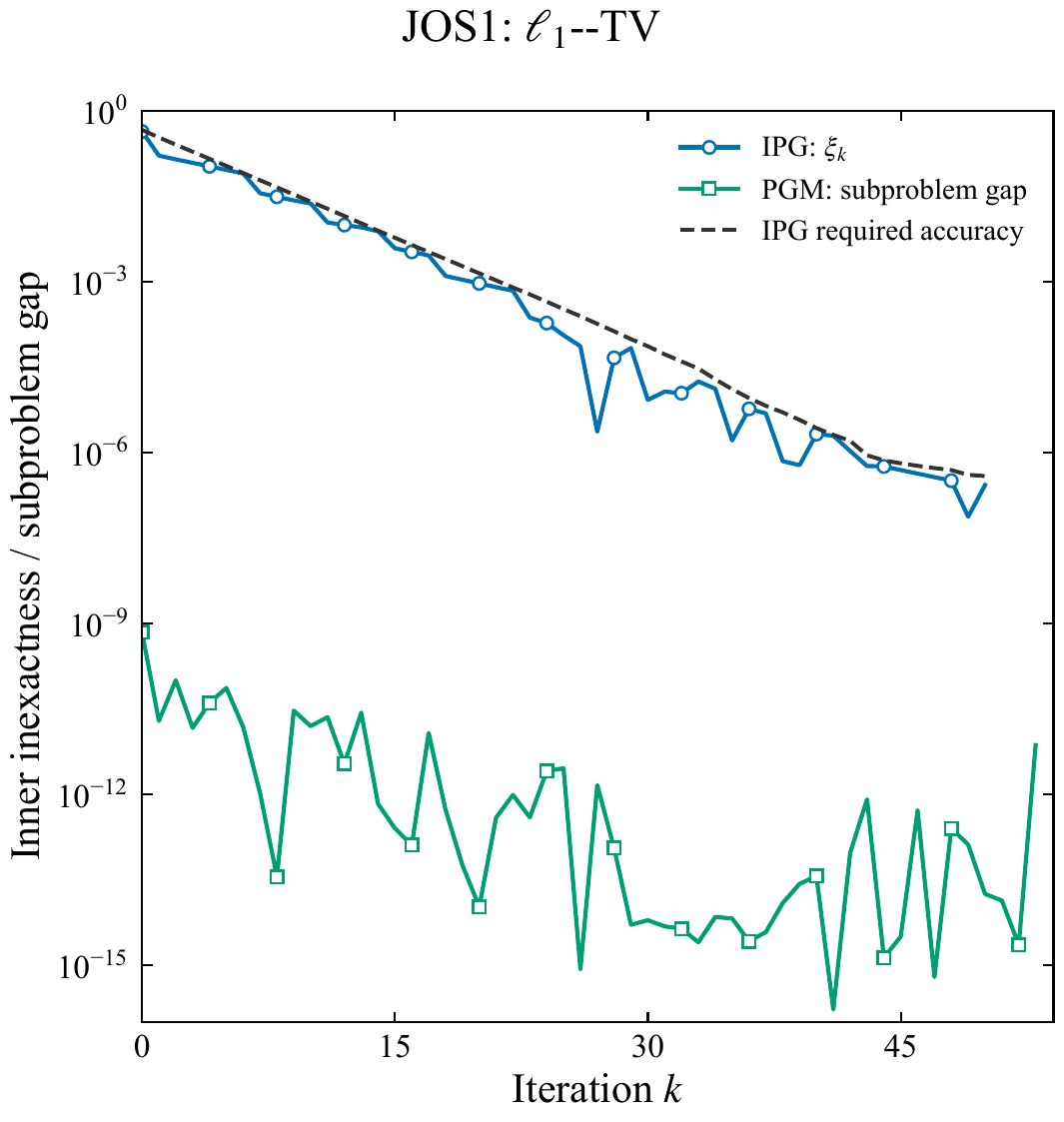}
	\caption{Inner errors for JOS1.}
	\label{fig:jos1-inner}
\end{minipage}\par
\addvspace{\intextsep}

\subsection{Regularized problems}

We next consider $F_i=f_i+g_i+h$, with
\begin{equation}
	\label{eq:regularized-model}
	g_i(x)=r_i^\top|Ax|,\qquad h(x)=\delta_{\mathcal C}(x),
	\qquad A=I+0.25S,
\end{equation}
where $S$ is the upper-shift matrix, $\mathcal C$ is a box, $\delta_{\mathcal C}$ denotes its indicator function, and $r_i$ are positive vectors. For simplex weights $\lambda\in\Delta_m$, the regularizer is
\[
\sum_{i=1}^m\lambda_i g_i(x)+h(x)
=(\sum_{i=1}^m\lambda_i r_i)^\top|Ax|
+\delta_{\mathcal C}(x),
\]
so IPG can exploit the common matrix $A$ in the calculation of the proximal
operator. The smooth objectives are based on JOS1~\citep{jin2001}, ZLT1-R~\citep{huband2006}, SLCDT2~\citep{schutze2008}, and CQ4-R, CQ8-R, and
CQ9-R~\citep{glasmachers2019}. Each CQ family has two and three objective variants. This gives nine configurations at each dimension. For each configuration, the two methods share ten initial points. We fix $\beta_k=5$ for CQ, $2$ for ZLT1-R, and $1$ for JOS1 and SLCDT2.

Table~\ref{tab:regularized-results} reports mean runtimes and outer iteration counts for nine problem settings at three dimensions. IPG achieves the lowest mean runtime in all $27$ configurations. At $n=500$, averaging the nine reported runtimes gives $0.6884$ seconds for IPG and $2.0693$ seconds for PGM, a reduction of approximately $67\%$. At $n=1000$, the corresponding averages are $1.5856$ and $6.0568$ seconds, respectively, giving a reduction of approximately $74\%$. The outer iteration counts are generally comparable. For example, for JOS1 at $n=1000$, IPG and PGM require
$606.5$ and $607.4$ outer iterations, respectively, whereas their mean runtimes are $0.6758$ and $15.4470$ seconds. These results are consistent with a lower computational cost per outer iteration for IPG.

\par\addvspace{\intextsep}
\noindent\begin{minipage}{\linewidth}
	\captionsetup{type=table}
	\centering
	\caption{Results for regularized test problems.}
	\label{tab:regularized-results}
	\footnotesize
	\setlength{\tabcolsep}{3pt}
	\renewcommand{\arraystretch}{0.93}
	\setlength{\defaultaddspace}{2pt}
	\begin{tabular}{@{}llrrrrrr@{}}
		\toprule
		& & \multicolumn{3}{c}{Mean runtime (s)} & \multicolumn{3}{c}{Mean outer iterations} \\
		\cmidrule(lr){3-5}\cmidrule(l){6-8}
		Problem ($m$) & Method & $n=300$ & $n=500$ & $n=1000$ & $n=300$ & $n=500$ & $n=1000$ \\
		\midrule
		JOS1 (2) & IPG & \textbf{0.3752} & \textbf{0.4961} & \textbf{0.7678} & 329.8 & \textbf{427.7} & \textbf{607.2} \\
		& PGM & 3.4729 & 6.4671 & 16.6854 & 329.8 & 427.8 & 607.4 \\
		\addlinespace
		ZLT1-R (5) & IPG & \textbf{0.1404} & \textbf{0.2404} & \textbf{0.8136} & \textbf{17.1} & \textbf{27.6} & \textbf{34.2} \\
		& PGM & 0.3168 & 0.8008 & 2.0447 & 19.4 & 32.3 & 35.9 \\
		\addlinespace
		SLCDT2 (3) & IPG & \textbf{1.6566} & \textbf{2.0216} & \textbf{3.1872} & 159.3 & 210.9 & 294.8 \\
		& PGM & 1.9353 & 3.7037 & 8.9054 & \textbf{156.6} & \textbf{205.0} & \textbf{276.6} \\
		\addlinespace
		CQ4-R (2) & IPG & \textbf{0.4456} & \textbf{0.8853} & \textbf{1.5518} & \textbf{72.3} & \textbf{97.8} & 135.0 \\
		& PGM & 0.8576 & 1.8308 & 4.8142 & 72.9 & 100.9 & \textbf{134.5} \\
		\addlinespace
		CQ4-R (3) & IPG & \textbf{0.5875} & \textbf{0.8556} & \textbf{2.2711} & \textbf{70.4} & \textbf{89.2} & \textbf{136.2} \\
		& PGM & 0.9589 & 1.9499 & 6.0630 & 72.5 & 97.3 & 140.4 \\
		\addlinespace
		CQ8-R (2) & IPG & \textbf{0.5345} & \textbf{0.7752} & \textbf{1.5156} & \textbf{70.2} & \textbf{91.1} & \textbf{130.0} \\
		& PGM & 0.9503 & 1.7225 & 4.8636 & 78.7 & 95.3 & 130.8 \\
		\addlinespace
		CQ8-R (3) & IPG & \textbf{0.5398} & \textbf{0.8281} & \textbf{1.9450} & \textbf{65.8} & \textbf{86.2} & \textbf{121.1} \\
		& PGM & 1.0097 & 1.8733 & 5.7713 & 74.8 & 92.8 & 129.1 \\
		\addlinespace
		CQ9-R (2) & IPG & \textbf{0.4275} & \textbf{0.4718} & \textbf{1.7517} & \textbf{66.5} & \textbf{90.3} & \textbf{121.7} \\
		& PGM & 0.8307 & 1.7159 & 4.4678 & 70.8 & 95.5 & 123.9 \\
		\addlinespace
		CQ9-R (3) & IPG & \textbf{0.5197} & \textbf{0.8456} & \textbf{2.0956} & \textbf{62.1} & \textbf{86.0} & \textbf{113.3} \\
		& PGM & 0.9005 & 1.9112 & 5.3968 & 67.3 & 96.4 & 122.7 \\
		\bottomrule
	\end{tabular}
\end{minipage}\par
\addvspace{\intextsep}

\subsection{Fairness-aware graphical estimation}

In this subsection, we consider a fairness-aware graphical estimation problem. Following~\citet{zhou2024fairness}, we seek a common matrix that keeps the overall loss small while reducing differences in group losses relative to fixed group specific references.

Let $q$ denote the number of groups and $p$ the number of variables. Let $S$ and $S_s$ be the sample covariance matrices of all samples and group $s$, respectively. The matrix variable belongs to $\mathcal C_p=\left\{\Omega\in\mathbb R^{p\times p}:\Omega=\Omega^\top,\ \Omega\succeq10^{-5}I\right\}$. For each group, let $\Omega_s^{\mathrm{loc}}$ be a fixed reference obtained by minimizing $L(\Omega,S_s)+\lambda\|\Omega\|_1$ over $\mathcal C_p$. Define the relative group loss and its squared differences from the other groups by
\[
e_s(\Omega)=L(\Omega,S_s)-L(\Omega_s^{\mathrm{loc}},S_s),\qquad D_s(\Omega)=\frac12\sum_{t=1}^{q}\bigl(e_s(\Omega)-e_t(\Omega)\bigr)^2.
\]
Thus, $D_s$ is small when the relative loss of group $s$ is close to those of the other groups. In particular, all $D_s$ vanish when $e_1=\cdots=e_q$.
We consider the following problem with $m=q+1$ objectives:
\[
\min_{\Omega\in\mathcal C_p}\bigl(F_1(\Omega),\ldots,F_{q+1}(\Omega)\bigr),
\]
where
\[
F_1(\Omega)=L(\Omega,S)+\lambda\|\Omega\|_1,\qquad F_{s+1}(\Omega)=D_s(\Omega)+\lambda\|\Omega\|_1,\quad s=1,\ldots,q.
\]
The first objective controls the loss on all samples, while the remaining objectives reduce differences between groups. The common regularization term promotes sparsity, with $\|\Omega\|_1=\sum_{i,j}|\Omega_{ij}|$ including diagonal entries. In~\eqref{def:MOP}, take $g_i(\Omega)=\lambda\|\Omega\|_1$ for every $i$ and $h=\delta_{\mathcal C_p}$.
We use two choices of the loss:
\[
L_{\mathrm{Cov}}(\Omega,S)=\frac12\|\Omega-S\|_F^2-0.01\log\det\Omega,\quad
L_{\mathrm{GL}}(\Omega,S)=-\log\det\Omega+\operatorname{tr}(S\Omega).
\]
Fair CovGraph estimates a covariance matrix, whereas Fair GLasso estimates its inverse. In both models, the terms common to $L(\Omega,S_s)$ and $L(\Omega,S_t)$ cancel, so $e_s(\Omega)-e_t(\Omega)$ is affine in $\Omega$. Consequently, each $D_s$ is a convex quadratic function.
The negative logarithm can be continued below the eigenvalue lower bound by a convex quadratic with matching value and derivative, giving an $L$-smooth extension of the losses without changing them on $\mathcal C_p$.

Let $\Omega^{\mathrm{std}}$ denote the standard reference obtained by minimizing $F_1$. The initial matrix is $\Omega^0=\Omega_{s_0}^{\mathrm{loc}}$. Group differences are summarized by $\Delta_{\mathrm{Cov}}=\frac1q\sum_{s=1}^{q}D_s$ and $\Delta_{\mathrm{GL}}=\sum_{s=1}^{q}D_s$.

Fair CovGraph uses the UCI Default of Credit Card Clients data with $p=20$, $q=8$ and $m=9$. We set $\lambda=0.1$ and $\beta_k=0.8$. Fair GLasso uses TCGA-LUAD data with $p=147$ genes and $q=2$ groups, giving $m=3$ objectives. We set $\lambda=0.03$, $\beta_k=0.5$, and a maximum of $100$ outer iterations. Table~\ref{tab:fair-results} reports $F_1$, runtime, outer iterations, and stopping residual $R_k$. 

\begin{table}[!htbp]
	\centering
	\begingroup
	\small
	\setlength{\tabcolsep}{6pt}
	\caption{Results of UCI Default of Credit Card Clients and TCGA-LUAD}
	\label{tab:fair-results}
	\begin{tabular}{@{}lrrrrr@{}}
		\toprule
		Method & $F_1$ & $\Delta$ & Time (s) & Outer & $R_k$ \\
		\midrule
		\multicolumn{6}{@{}l}{Fair CovGraph ($p=20$, $q=8$, $m=9$)} \\
		Standard & 9.405913 & 0.530824 & -- & -- & -- \\
		IPG & 9.438920 & 0.283925 & \textbf{0.0629} & 11 & $2.43\times10^{-5}$ \\
		PGM & 9.430211 & 0.304471 & 0.3339 & \textbf{7} & $4.82\times10^{-5}$ \\
		\midrule
		\multicolumn{6}{@{}l}{Fair GLasso ($p=147$, $q=2$, $m=3$)} \\
		Standard & 74.896979 & 1.746164 & -- & -- & -- \\
		IPG & 76.155449 & 0.002171 & \textbf{89.0223} & 100 & 0.504805 \\
		PGM & 75.439404 & 0.007163 & 531.1115 & 100 & 0.818466 \\
		\bottomrule
	\end{tabular}
	\endgroup
\end{table}

For Fair CovGraph, IPG and PGM reduce $\Delta$ relative to Standard by approximately $43.7\%$ and $42.6\%$, respectively, while increasing $F_1$ by only $0.28\%$ and $0.26\%$. Both stopping residuals are below $10^{-4}$. IPG requires $16$ outer iterations, compared with $7$ for PGM, with recorded runtimes of $0.1545$ and $0.4497$ seconds, respectively. For Fair GLasso, both methods reduce $\Delta$ by more than $99.6\%$ after $100$ outer iterations, with increases in $F_1$ of approximately $0.34$--$0.35\%$. IPG attains smaller values of $F_1$, $\Delta$, and $R_k$ than PGM at this iteration limit.

\section{Conclusion}
\label{sec:conclusion}
We proposed an inexact proximal gradient method for composite multiobjective optimization based on an approximate optimality condition and a hybrid stopping criterion. The subproblems are solved using an inexact projected gradient method and ADMM based on the proximal operators of $g_i$ and $h$. Under the stated assumptions, we established convergence of the objective values and weak Pareto optimality of every accumulation point, together with a convergence rate $\mathcal O(1/k)$  under weighted summability of the errors. The numerical results retained from the original implementation illustrate the computational benefits of inexact subproblem solutions.

Future work will focus on Nesterov acceleration, including stopping criteria and inner stopping rules that support accelerated convergence while controlling the cost of inexact computations.

\section*{Acknowledgements}

\bibliographystyle{apalike-ejor}
\bibliography{sample-ejor}

\end{document}